\documentclass[a4paper]{siamart250211}

\usepackage{damacros}
\usepackage{subcaption}

\usepackage{bm}

\usepackage{cite}

\usepackage{tikz} 
\usepackage{caption}

\renewcommand{\vec}[1]{\bm{#1}}

\newcommand{\vecc}{\vec{c}}
\newcommand{\vecd}{\vec{d}}

\newcommand{\vecf}{\vec{f}}
\newcommand{\vecw}{\vec{w}}
\newcommand{\vecx}{\vec{x}}
\newcommand{\vecy}{\vec{y}}

\newcommand{\LL}{\mathcal{L}}
\newcommand{\RR}{\mathbb{R}}

\newsiamremark{remark}{Remark}
\crefname{remark}{Remark}{Remarks}
\Crefname{remark}{Remark}{Remarks}
\crefname{hypothesis}{Hypothesis}{Hypotheses}
\Crefname{hypothesis}{Hypothesis}{Hypotheses}

\newcommand{\rev}[1]{{#1}} 

\headers{RBFs for neural field models}{Sage Shaw, Zachary P Kilpatrick, Daniele Avitabile}

\title{Radial Basis Function Techniques for Neural Field Models on Surfaces} 

\author{
  Sage B. Shaw\thanks{Department of Applied Mathematics, University of Colorado Boulder, Boulder, CO 80309 USA
    (\email{Sage.Shaw@colorado.edu})}
  \and 
  Zachary P. Kilpatrick\thanks{Department of Applied Mathematics, University of Colorado Boulder, Boulder 80309 CO, USA
    (\email{zpkilpat@colorado.edu})}
  \and
  Daniele Avitabile\thanks{Department of Mathematics, Vrije Universiteit Amsterdam,
    The Netherlands, and Inria MathNeuro team, Montpellier, France
    (\email{d.avitabile@vu.nl}) \\[1ex] {\bf Funding:} This work was initiated through discussions at the Institute for Computational 
and Experimental Research in Mathematics (ICERM), which is supported by the National 
Science Foundation under grant DMS-1929284.
S.B.S. \& Z.P.K. were supported by the National Science Foundation under grant 
DMS-2207700. 
}
}

\graphicspath{{./Figures/}}

\begin{document}

\maketitle

\begin{abstract}
  We present a numerical framework for solving neural field equations on surfaces
  using Radial Basis Function (RBF) interpolation and quadrature. Neural field models
  describe the evolution of macroscopic brain activity, but modeling studies often
  overlook the complex geometry of curved cortical domains. Traditional numerical
  methods, such as finite element or spectral methods, can be computationally
  expensive and challenging to implement on irregular domains. In contrast, RBF-based
  methods provide a flexible alternative by offering interpolation and quadrature
  schemes that efficiently handle arbitrary geometries with high order accuracy. We
  first develop an RBF-based interpolatory projection framework for neural field
  models on general smooth surfaces.
  Quadrature for both flat and curved domains are derived in detail, ensuring
  high-order accuracy and stability as they depend on RBF hyperparameters (basis
  functions, augmenting polynomials, and stencil size). Through numerical
  experiments, we demonstrate the convergence of our method, highlighting its
  advantages over traditional approaches in terms of flexibility and accuracy. We
  conclude with an exposition of numerical simulations of spatiotemporal activity on
  complex surfaces, illustrating the method’s ability to capture complex wave
  propagation patterns.
\end{abstract}

\section{Introduction}\label{sec:introduction} 

Neural field equations are nonlinear integro-differential equations that model large-scale neuronal activity, offering tractable formulations for analysis and simulation across diverse domains~\cite{bressloff2011spatiotemporal}. Such activity is strongly determined by the architecture of synaptic connections between neurons~\cite{coombes2010large}. Associated models typically grapple with the high volume of neurons and synapses, especially in primate cortex, with approximations such as coarsening, homogenization, and timescale separations~\cite{breakspear2017dynamic}. Neural fields treat cortex as a continuum excitable medium and model connectivity via a spatial weight kernel in a nonlinear integral operator, allowing the application of theory from nonlinear waves and other methods from partial differential equations (PDE)~\cite{bressloff2011spatiotemporal}. The kernel defines a spatially dependent synaptic coupling strength based on the locations of pre- and post-synaptic neurons. Analyses and simulations typically center on this nonlocal term, which drives local dynamics throughout the domain.


Standard derivations of neural fields are grounded in the known spatial and functional structure of synaptic connectivity throughout cerebral cortex, justifying the typical spatial coarse graining required to obtain continuum models~\cite{wilson1973mathematical,shipp2007structure}. For analytical and numerical convenience, most modeling studies consider canonical spatial domains and connectivity functions (e.g., distance-dependent weight kernels on planar domains)~\cite{coombes2014neural}. Amari (1977)~\cite{amari1977dynamics} introduced a simplified Heaviside transfer function to enable the use of interface methods and explicit construction of wave solutions. Weight kernels with rational Fourier transforms reduce nonlocal neural fields to local PDEs, facilitating efficient simulation and analysis~\cite{laing2003pde}. However, these simplifications often ignore geometrical complexity. Here we advance a method that is highly flexible to the formulation of neural fields on arbitrary domains with curvature and nonstandard weight kernels. 


Although cortex is three-dimensional, many imaging modalities (e.g., fMRI, VSD) access activity on its two-dimensional surface~\cite{lee2005traveling,huang2010spiral}. Prior models considering curved domains appeal to idealizations like spheres~\cite{bressloff2003spherical,visser2017standing} and tori~\cite{kneer2014nucleation}, allowing for spectral solution methods leveraging known orthonormal eigenbases. General curved surfaces that cannot be expressed as Cartesian products lack this structure. Recent work has explored neural fields on arbitrary smooth surfaces~\cite{martin2018numerical}, but without convergence guarantees and with only first-order accuracy. A recent framework for projection-based neural field solvers~\cite{Avitabile:2023ab} offers an operator-theoretic approach to convergence analysis, decomposing numerical error into contributions from spatial projection, time integration, and quadrature. This unified perspective applies to both Galerkin and collocation methods, clarifying how errors scale with resolution and guiding the design of efficient discretizations for simulations on complex geometries.

In this work, we develop a high-order numerical method for simulating neural field equations on complex geometries and advance the associated convergence theory within a projection-based framework. Our approach combines radial basis function quadrature (RBF-QF)~\cite{REEGER2016,REEGER2018,REEGER2016Sphere} with RBF interpolation to achieve high-order accuracy using relatively few degrees of freedom. This is especially important for simulations on curved surfaces, where evaluating nonlocal integral operators at $\mathcal{O}(n^2)$ cost can be prohibitive.\footnote{On certain domains, FFT-based acceleration reduces complexity to $\mathcal{O}(n \log n)$, but such approaches are not applicable to general manifolds.} We formulate projection and collocation schemes for neural fields and present convergence results for projection methods (\Cref{sec:model}); we describe the RBF-QF algorithm and its extension to manifolds (\Cref{sec:rbfs}); and we verify performance through numerical tests and simulations on nontrivial geometries (\Cref{sec:numerical_tests,sec:showcase}).

\section{Collocation schemes for the neural field model}\label{sec:model} 
We consider, as a model problem, a neural field posed on a compact domain $(t,\vecx)
\in [0,T] \times (D \subset \RSet^{\textrm{dim}})$:
\begin{equation} \label{eq:NF}
    \begin{aligned}
      \partial_t u(t, \vecx) &= - u(t,\vecx) + g(t,\vecx) + \int_D w(\vecx,\vecy) f(u(\vecy,t)) d \mu (\vecy), \\
        u(0,\vecx) &= v(\vecx).
    \end{aligned}
\end{equation}
We assume the integral is on a volume or surface $D$, and so is expressed in terms
of a measure $\mu$. Heuristically, a collocation scheme for \cref{eq:NF} can be
conceived in two steps, as follows: \\
\vspace{-4mm}

\noindent
{\bf Step 1: collocation.} Select $n$ points $\Xi = \{ \vecx_i : i
    \in \NSet_n := \{1,...,n\} \} \subset D$ in the domain, where $n$ {\em may} depend on a discretization parameter $h$ so $n(h) \to \infty$ as $h \to 0$.
  We impose that \cref{eq:NF} holds at each node, that is,
  \textit{we collocate the equation} on $\Xi$. This leads to a functional equation
  relating $\{ \partial_t u(t,\vecx_i) \}_i$ and $\{ u(t,\vecx_i)\}_{i}$ to the
  function $u(t,\blank)$, appearing under the integral of \cref{eq:NF}.  \\
  \vspace{-3mm}
  
\noindent
{\bf Step 2: quadrature.} Select a quadrature for the integral, in terms of
    $\{ u(t,\vecx_i) \}_i$ and obtain a system of coupled, nonlinear ordinary
    differential equations (ODEs). \\
\vspace{-3mm}

Avitabile (2023) \cite{Avitabile:2023ab} builds on previous work by
Atkinson (2005)~\cite{atkinson2005theoretical} to show that convergence estimates
for collocation (and other) schemes are still possible when disregarding {\bf Step
2}. {\bf Step 1} is a
fully functional scheme for which one can estimate convergence rates of the
collocation (or other projection types) as $n \to \infty$. {\bf Step 2} further
approximates the projection with quadrature, establishing a new scheme (a
\textit{discrete collocation scheme} in Atkinson's
terminology~\cite{atkinson2005theoretical}).

\rev{
Henceforth we assume the following hypotheses,  provide a natural functional
setup for the neural field problem posed on $\XSet = C(D)$ \cite{Potthast:2010kb, Avitabile:2023ab}
\begin{hypothesis}[General hypotheses]\label{hyp:general}
  \begin{enumerate}
    \item \label{hyp:domain}
      The cortical domain $D \subset \RSet^{\textrm{dim}}$ is compact. 
  \item \label{hyp:timeInterval}
    The temporal domain $J \equiv [0,T] \subset \RSet$ is compact. 
  \item \label{hyp:kernel}
    The synaptic kernel $w \colon D \times D \to \RSet$ is a function in $C(D \times
    D)$.
  \item The firing rate $f \colon \RSet \to \RSet$ is a bounded and
  everywhere differentiable Lipschitz function, which guarantees $f,f' \in B(\RSet)$.
  \item \label{hyp:initialCondition}
    The initial condition $v \colon D \to \RSet$ is a function in $\XSet = C(D)$.
  \item The forcing $g \colon D \times J \to \RSet$ is a function in
    $C(J,\XSet)$.
  \end{enumerate}
\end{hypothesis}

To set the stage for our collocation scheme, we rewrite the neural field \cref{eq:NF} in operator form
\begin{equation}\label{eq:cauchy}
  \begin{aligned}
  & u'(t) = - u(t) + WF(u(t)) + g(t) =:  N(t,u(t)), \qquad t \in J \equiv [0,T], \\
  & u(0) = v,
  \end{aligned}
\end{equation}
where we define the integral operator $W$ and nonlinear operator $F$:
\begin{equation} \label{eq:WFDef}
  \begin{aligned}
   W \colon \XSet \to \XSet,
   \quad 
   v \mapsto \int_{D} w(\blank,\vec{y}) v(\vec{y}) \,d \mu(\vec{y}),
   \qquad 
   F \colon \XSet \to \XSet, 
   \quad 
   v \mapsto f(v).
  \end{aligned}
\end{equation}
Note, \cref{eq:NF} evolves the $\RSet$-valued function $u
\colon D \times J \to \RSet$, whereas \cref{eq:cauchy} evolves
the $\XSet$-valued function $U \colon J \to \XSet$, $t \mapsto u(\blank,t)$. With
a small abuse of notation, we adopt the same symbol for both functions ($u$), as for
the forcing function $g$.

\Cref{eq:cauchy} describes a Cauchy problem on the Banach space $\XSet$, which is well-posed under
\cref{hyp:general}. That is, there exists a unique $u
\in C^1(J,\XSet)$ satisfying \cref{eq:cauchy}, as shown in \cite[Lemma
2.7]{Avitabile:2023ab} (See \cite{Potthast:2010kb} for analogous results on unbounded
cortices). 
}

\rev{
\subsection{Step 1: collocation via RBF interpolatory projection}\label{ssec:projectionSchemes} 
RBFs provide a flexible and accurate way to interpolate functions and have been
successfully applied to nonlinear PDEs with complex spatiotemporal dynamics, like
reaction–diffusion systems~\cite{shankar2015radial}. Since our goal is to approximate
nonlinear integrodifferential (neural field) equations, RBF-based collocation schemes
are a natural choice. We present these schemes using the abstract notion of an
\textit{interpolating projector}, which maps functions onto their interpolants
(e.g.,~RBFs).

In this setup, one seeks a solution $u$ to the neural field equations as a mapping
on $[0,T]$ to $\XSet = C(D)$, the space of continuous functions on $D$. An
interpolatory projector $P_n$ linearly projects any function $v \in \XSet$ to
a function $v_n = P_n v$ in an $n$-dimensional subspace $\XSet_n = \spn\{L_1, \ldots, L_n\}
\subset \XSet$, with the property
\begin{equation}\label{eq:PnInterp}
  P_n \colon \XSet \to \XSet_n, \qquad P_n v = v, \qquad \text{on $\Xi$.}
\end{equation}
The approximant then takes the form
\[
  v_n(\vecx) := (P_n v)(\vecx) = \sum_{j=1}^{n} v(\vecx_j) L_j(\vecx), \qquad x \in D.
\]
}

\rev{
The collocation scheme approximates the solution $u$ of \cref{eq:cauchy}, by a
function $u_n \in C^1(J,\XSet_n)$, that satisfies the following
Cauchy problem on the $n$-dimensional subspace $\XSet_n \subset \XSet$:
\begin{equation}\label{eq:cauchyPn}
  \begin{aligned}
  & u'_n(t) = - u_n(t) + P_n WF(u_n(t)) + P_n g(t) =  P_n N(t,u_n(t)), \qquad t \in J, \\
  & u_n(0) = P_n v,
  \end{aligned}
\end{equation}
that is, the Cauchy problem obtained from \cref{eq:cauchy} upon applying the projection $P_n$ to the
vectorfield $N$, and to the initial condition $v$. 

In view of \cite[Proposition 4.1]{Avitabile:2023ab}, this is equivalent to imposing that \cref{eq:NF} holds
at all $\vecx_i \in \Xi$ (hence we have collocated the equation, and completed {\bf Step 1}),
\[
  \begin{aligned}
  & u'_n(t)(\vecx_i) = \bigl[P_n N(t,u_n(t))\bigr](\vecx_i), \qquad (i,t) \in \NSet_n \times J, \\
  & u_n(0)(\vecx_i) = (P_n v)(\vecx_i).
  \end{aligned}
\]
In \cite[Theorem 3.1]{Avitabile:2023ab} it is shown that, under \cref{hyp:general},
also the projected problem is well-posed, admitting a unique solution $u_n \in C^1(J,\XSet_n)$.
}

\rev{
The projected provlems \cref{eq:cauchyPn} evolves the function $u_{n}(t)$ in
$\XSet_n$. This form is useful in the analysis of the scheme, but for
implementations we pursue an ODE approximating the $n$ coefficients $\alpha_i (t)
:= u_n(t)(\vec{x}_i)$. 
%
Since $u_n(t) \in \XSet_n$, we set $u_n(t) = \sum_{j \in \NSet_n}
\alpha_j(t) L_j$, exploit the Lagrange property of the basis $\{ L_j \}$, and obtain an ODE
in $\RSet^n$
\begin{equation} \label{eq:NFODE}
  \begin{aligned}
   & \alpha'(t) = -\alpha(t) + K(\alpha(t)) + \gamma(t) \qquad  t \in J,\\ 
   & \alpha (0) = \nu
  \end{aligned}
\end{equation}  
where 
\[
  \begin{aligned}
  & \alpha_i(t) = u_n(t)(\vecx_i), && \nu_i = v(\vecx_i),  \\
  & \gamma_i(t) = g(\vecx_i,t),
  && K_i(\alpha) = 
  \int_{D}w(\vecx_i,\vec{y}) 
  f\Bigl( \sum_{j \in \NSet_n} \alpha_j L_j(\vec{y})\Bigr) \, d\mu(\vec{y}).
  \end{aligned}
\]
We stress that \cref{eq:cauchyPn} and \cref{eq:NFODE} are equivalent Cauchy problems
on the $n$-dimensional phase spaces $\XSet_n$ and $\RSet^n$, respectively.
}

\subsection{Step 2: quadrature and Discrete Projection Schemes} 
The 

\noindent The ODE \cref{eq:NFODE} on $\RSet^n$ is not yet implementable on a computer, because
generally the integrals in $K$ can not be evaluated in closed form, and must be
approximated by numerical quadrature. Such an approach generates a computationally implementable algorithm using \textit{Discrete Projection
Schemes}.

We use a quadrature scheme induced by interpolation. We will later discuss the
details of the quadrature, which approximates the integral operator
\[
  Q \colon \XSet \to \RSet,
  \qquad 
  Qv = \int_{D} v(\vec{y}) \,d\mu(\vec{y}) 
\]
with a sum of the form
\[
  Q_n \colon \XSet \to \RSet,
  \qquad 
  Q_n v = \sum_{j=1}^{n} v(\vec{x}_j) \mu_j.
\]

We can apply the quadrature scheme $Q_n$ to the operator $W$ in \cref{eq:cauchyPn},
define
\begin{equation}\label{eq:WnDef}
  W_n \colon \XSet \to \XSet_n, \qquad (W_n v)(\vec{x}) = Q_n (
  w(\vec{x},\blank) v),
\end{equation}
and arrive at the scheme
\begin{equation}\label{eq:cauchyPnQuad}
  \begin{aligned}
  & \tilde u'_n(t) = - \tilde u_n(t) + P_n W_nF(\tilde u_n(t)) + P_n g(t), \qquad t \in J, \\
  & \tilde u_n(0) = P_n v,
  \end{aligned}
\end{equation}
the ODE system for $u$ at the collocated points, approximated by quadrature.

Following similar steps to \Cref{ssec:projectionSchemes}, 
we 
obtain
\begin{equation} \label{eq:NFODETilde}
  \begin{aligned}
   & \tilde \alpha'(t) = -\tilde \alpha(t) + \tilde K(\tilde \alpha(t)) + \gamma(t) \qquad  t \in J,\\ 
   & \alpha (0) = \nu,
  \end{aligned}
\end{equation}  
for the evolution of the vector of coefficients $\alpha(t) = [\alpha_1(t), ..., \alpha_n(t)]^T$, where 
\[
  \begin{aligned}
  & \tilde \alpha_i(t) = \tilde u_n(t)(\vecx_i), && \nu_i = v(\vecx_i),  \\
  & \gamma_i(t) = g(t,\vecx_i),
  && \tilde K_i(\tilde \alpha) = 
    \sum_{j = 1}^n w(\vecx_i,\vecx_j) f\Bigl( \sum_{k \in \NSet_n} \tilde \alpha_k
    L_k(\vecx_j)\Bigr)
    \mu_j.
  \end{aligned}
\]

\subsection{Convergence of the collocation scheme}\label{ssec:convEstimate} 
In \cite{Avitabile:2023ab}, Atkinson's approach for studying projection methods in
Fredholm and Hammerstein integral equations
\cite{atkinson1997,atkinson2005theoretical}, is extended to time-dependent neural
field equations. The central result of the paper is a bound for the error $\| u -
u_{n} \|_{C(J,\XSet)}$ in terms of the projection error $\| u - P_n
u\|_{C(J,\XSet)}$, that is, a bound for the collocation scheme in {\bf Step 1}, in terms of
the error of the interpolating projector. 

\begin{theorem}[Abridged from \cite{Avitabile:2023ab}, Theorem 3.3 and discussion on
  page 570]\label{thm:convergence}
Under \cref{hyp:general}, if $P_n v \to v$ for all $v \in \XSet$, then $u_n \to u$ in
$C(J,\XSet)$. Further, there exist positive constants $m$ and $M$, independent of $n$, such that
\[
  m\| u - P_n u  \|_{C(J,\XSet)} \leq \| u - u_n  \|_{C(J,\XSet)} \leq M \| u - P_n u  \|_{C(J,\XSet)}.
\]
\end{theorem}

Since the time interval $J$ is compact it holds $\| u - P_n u \|_{C(J,\XSet)} = \|
u(t_*) - P_n u(t_*)\|_{\XSet}$ for some $t_* \in J$, and thus the convergence rate of the
collocation scheme is estimated from the interpolation error $\| v - P_n v
\|_{\XSet}$ for some $v \in \XSet$.

The error of the Discrete Collocation Scheme can in principle be estimated, given
$\| u - \tilde u_n \| \leq \| u - u_n \| + \| u_n - \tilde u_n \|$, but we
do not take this route here. Controlling $\| u_n - \tilde u_n \|$ requires rigorous
error estimates for quadrature rules on curved domains $D$ which, to the best
of our knowledge, are not yet available when the underlying interpolant is an RBF.
Instead, we proceed in the remaining sections as follows:
\begin{enumerate}
  \item We employ RBF interpolants whose interpolation error $\| P_n v - v \|$ is
    available in literature, and deduce rates for the error of collocation scheme
    $\| u - u_n \|$, as per \cref{thm:convergence}.
  \item We select a quadrature scheme built on RBFs, for which we provide numerical
    evidence of convergence rates. A guiding principle is to select a quadrature rule that
    matches the error of the collocation scheme.
  \item We provide numerical evidence of convergence rates for the error $\| u - \tilde u_n \|$ of
    the fully discrete scheme in cases where an analytical solution to the neural
    field is available in closed form.
\end{enumerate}

\section{RBF interpolation and quadrature} \label{sec:rbfs}
We use RBFs for the interpolatory projector $P_n$ and quadrature formula $Q_n$
introduced in \Cref{sec:model}. In particular, we use RBF quadrature formulae
(RBF-QF), which provide an interpolation-based technique to calculate quadrature
weights for arbitrary sets of quadrature nodes with high-order accuracy.
This quadrature is geometrically flexible: it works in domains of arbitrary dimension and
in complex geometry. The approximation power comes from the remarkable properties
of radial basis function (RBF) interpolation. In this section, we describe the
algorithm used to interpolate functions, generate quadrature weights in flat and
curved domains, and discuss the consequences of hyperparameter choices. This general approach can then be leveraged to implement our neural field approximation scheme.

\subsection{Local Interpolation} \label{ssec:local_interpolation}
The first use of RBF interpolation was by Hardy in 1971~\cite{HARDY1971} to construct
a topography from scattered elevation measurements. The interpolant was a linear
combination of radially-symmetric basis functions called {\em Multiquadrics}. Critically,
the centers of these basis functions coincided with the locations of the
measurements, and as a result, the interpolants were guaranteed to exist and be
unique. The Mairhuber-Curtis theorem shows that such guarantees break down when basis
functions are chosen independently of the interpolation nodes~\cite{FASSHAUER2007,
MAIRHUBER1956, CURTIS1959} as occurs in multivariate polynomial
interpolation. Later developments showed that other radially symmetric basis functions
could be used similarly, with resulting interpolants exhibiting spectral
convergence as the number of nodes increased ~\cite{buhmann1993spectral,
yoon2001spectral}. For a comprehensive introduction to RBF interpolation, see \cite{FASSHAUER2007}.

The price of spectral convergence, however, is that finding interpolant coefficients
requires solving an $n \times n$ linear system, which becomes prohibitively expensive and
ill-conditioned as $n$ (the number of interpolation nodes) increases. There are two
approaches to remedy this limitation. First, one can choose compactly supported basis
functions, often called
Wendland functions~\cite{FASSHAUER2007, WENDLAND2005}, to produce a sparse, banded linear
systems, which can be solved efficiently. We opt for the second
approach: local RBF interpolation using polyharmonic splines and appended polynomial
basis functions, which is the focus of this subsection. 
\rev{
For each stencil, unisolvency of the appended polynomial basis is 
ensured by selecting distinct, well-distributed nodes, so that the associated 
Vandermonde matrix is full rank. On curved surfaces, we work in local coordinate 
parameterizations of each stencil, so the same unisolvency conditions as in the 
planar case apply~\cite{bayona2024meshfree}.
}

\rev{
RBF interpolation is often touted as a ``mesh-free'' method and this is certainly
true of global RBF interpolation, which yields a smooth interpolant over the entire
domain. Local RBF interpolation is arguably mesh-free as well, however it produces a
piecewise interpolant that is smooth over each element in a partition of the domain
but potentially discontinuous along the boundaries of these elements. 
Generally, this partition is taken to be the Voronoi partition and is not explicitly calculated. 
It is thus mesh-free in the sense that the mesh is implicit. We emphasize this
because the RBF-QF algorithm for surfaces will explicitly require a triangulation in
lieu of a Voronoi partition, and is thus not a mesh-free algorithm. We now provide a
self-contained description of the interpolation algorithm.}

Given a compact domain $D \subset \mathbb{R}^{\textrm{dim}}$, we seek an approximation to a
sufficiently smooth function $f: D \to \mathbb{R}$, using a discrete finite set of interpolation nodes $\Xi_n \subset D$, and 
a partition of the domain $\{E_j\}_{j=1}^m$. We will approximate $f$ piecewise on
each element $E_j$, where each piecewise component will be a linear combination of
basis functions. For each element $E_j$, there is an associated subset of
interpolation nodes $S(E_j) \subseteq \Xi_n$, which we call a stencil. Generally, the
elements $E_j$ are small polygons that shrink or are divided as more points are
added. The stencil $S(E_j)$ contains the closest $k$ points in $\Xi_n$ to a
suitably-defined center of $E_j$, making $k = |S(E_j)|$ a hyper-parameter called the
{\em stencil size}. In what follows, we will sometimes use $E$ to refer
to a generic element in the partition $\{ E_j \}_{j =1}^m$.

\begin{figure}
  \centering
  \includegraphics[width=1.0\textwidth]{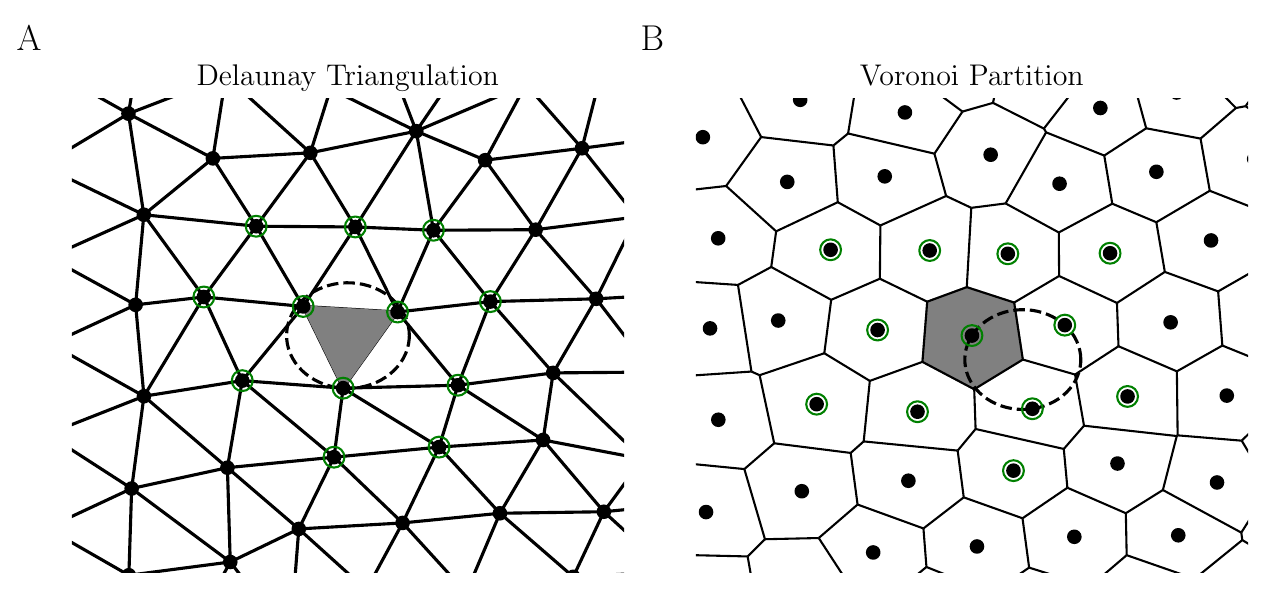}
  \caption{
  \textbf{Stencils.} Each panel shows a subset of the interpolation nodes $\Xi_n$
  (black dots) for a domain $D \subset \RSet^2$, the boundaries (solid lines)
  enclosing partition elements $E_j$, and an example element (shaded gray) and its
  stencil points (green circles). \textbf{A.}~Delaunay triangulations partition
  elements using interpolation nodes as vertices. An element's circumscribing circle
  does not
  contain any additional interpolation nodes. \textbf{B.}~Voronoi
  partitions assign one interpolation node per element. A circle centered at a
  Voronoi vertex can always be drawn to contain only the centers of neighboring
  elements.
  }
  \label{fig:stencil}
\end{figure}

Partitions are demonstrated in two examples depicted in \Cref{fig:stencil}. The Delaunay triangulation is popular for partitioning planar domains and widely integrated into many standard mathematical software libraries~\cite{fortune2017voronoi}. Interpolation nodes are used as the vertices of triangular elements that form the partition (\Cref{fig:stencil}{A}). The defining feature of the Delaunay triangulation is that a circle circumscribing each triangular element contains no other interpolation nodes except its vertices. Alternatively, Voronoi partitions (\Cref{fig:stencil}{B}) assign each element a single interpolation node and define it as the set of points in $D$ closer to that node than to any other. 
The Voronoi diagram is the graph dual of the Delaunay triangulation for a given set of points\footnote{This is almost always true. Non-uniqueness can arise when a Voronoi vertex is shared by more than three cells (e.g., in a Cartesian grid), but this is rare for randomly chosen points and can be resolved by small perturbations.}.

We will refer to $\{ E_j \}_{j=1}^m$ as the partition, mesh, or triangulation interchangeably
and
the subsets $E_j$ as elements (of the partition), patches, or triangles depending on
the context. While our presentation allows for general choices of interpolation nodes and meshes, not all such combinations yield accurate or stable interpolants.
A useful quantity associated with the nodes is the \textit{covering radius}~\cite{FASSHAUER2007,
buhmann2000} 
\begin{equation} 
  h = \sup_{\vecx \in D} \min_{\vecy \in \Xi_n} \|\vecx - \vecy \|,\label{eqn:spacing} 
\end{equation} 
the supremum of the radii of circles with centers $\vecx \in D$ that contain no
interpolation nodes. The covering radius is analogous to mesh spacing and scales asymptotically as $h \sim n^{-1/\text{dim}}$, as in regular grids. 

This notion arises naturally in both the Delaunay triangulation
and Voronoi partition. In the Delaunay case, the radius of any element's circumscribed circle
(See \Cref{fig:stencil}{A}) provides a lower bound on $h$, and if the triangulation covers $D$, then $h$ is the maximum of these radii across all triangular elements. In the Voronoi case, the
covering radius is the maximum distance from a node to a vertex of its associated patch. \Cref{fig:stencil}{B} shows a circle
centered on a Voronoi vertex with radius equal to the distance to the adjacent
nodes. In both cases, the covering radius is the maximum over a finite
set of such radii, making it computationally straightforward.

Going forward, we assume each stencil has fixed size $|S(E)| = k$ for
all $E \in \{E_j\}_{j=1}^m$. 
We introduce a parameter, ``$\textrm{deg}$'', to control the order of
accuracy of our interpolant. Over each element $E$, the
approximation is expressed as a linear combination of RBFs $\phi_{\vecy}$ and
polynomial basis functions $\pi_{\vec{\alpha}}(\vecx) = \vecx^{\vec{\alpha}}$, where
$\vec{\alpha}$ is a
multi-index, and $\{\pi_{\vec{\alpha}}\}_{|\vec{\alpha}| \le \deg}$ forms a basis for
$\mathbb{P}_{\text{deg}}(D)$, the space polynomials of degree at most $\text{deg}$.
A function $f \in C(D)$ is then approximated using the piecewise-defined function
\begin{equation}\label{eq:sFunDef}
  s(\vecx): =  \sum_{j = 1}^{m} s_{E_j}(\vecx) \chi_{E_j}(\vecx), 
  \qquad \vecx \in D,
\end{equation}
where $\chi_E \colon D \to \RSet$ is the indicator function on $E$ 
and where $s_E$ is the function
\[
  s_E(\vecx) = \sum_{\vecy \in S(E)} c_{E,\vecy} \phi_{\vecy}(\vecx) 
  + \sum_{|\vec{\alpha}| \le \text{deg}} d_{E, \vec{\alpha}} \pi_{\vec{\alpha}}(\vecx).
\]
The RBFs $\psi_{\vec{y}}$ are each translates of the same \textit{basic function} $\psi_{\vec{y}}(\vec{x}) = \Phi(\|\vecx - \vecy\|)$~\cite{FASSHAUER2007}, but
not all radially symmetric functions are suitable basic functions. We will restrict
ourselves to polyharmonic splines with fixed order parameter $\ell$
\begin{equation}\label{eq:polyHarmonic}
  \Phi(r) = 
  \begin{cases}
    r^\ell, & \text{for $\ell$ odd}, \\
    r^\ell \log(r), & \text{for $\ell$ even.}
  \end{cases}
\end{equation}

To determine the coefficients $c_{\vecy}$ and $d_{\vec{\alpha}}$ we enforce two sets of conditions known as
\textit{interpolation conditions} and \textit{moment conditions}, respectively:
\begin{align}
  f(\vecx) &= s_E(\vecx), 
  && \text{for $\vecx \in S(E)$,} 
    \tag{interpolation} \\
  0 & = \sum_{\vecy \in S(E)} d_{\alpha} \vecy^\alpha, 
  && \text{for $|\vec{\alpha}| \le \text{deg}$.}
    \tag{moment}
\end{align}
Note that though we refer to these as ``interpolation conditions'', a point $\vecx \in
S(E)$ but $ \vecx \not \in E$ is interpolated by the local interpolant $s_E$, but
not by the global piecewise interpolant $s$. The full approximation
interpolates $f$ at each $\vecx \in \Xi$ since every node lies in some element
$E$, but each local interpolant $s_E$ satisfies additional conditions not directly reflected in the global approximation. 

The moment conditions ensure \textit{polynomial reproduction}: for any $p \in \mathbb{P}_{\text{deg}}(D)$, the local interpolant satisfies $s_E = p$ on its domain. This property underpins the approximation power of the local interpolation scheme.

Thus, the coefficients are determined by solving the linear system
\begin{align}
    \begin{bmatrix}
        A & \Pi \\
        \Pi^T & 0 \\
    \end{bmatrix}
    \begin{bmatrix}
        \vecc \\ \vecd 
    \end{bmatrix}
    =
    \begin{bmatrix}
        \vecf_E \\
        \vec{0}
    \end{bmatrix} \label{eqn:interp_system}
\end{align}
where the full block matrix is the \textit{interpolation matrix}. Here, $A_{ij} = \Phi(\|\vecx_i - \vecx_j\|)$ is the \textit{RBF interpolation matrix}, $\Pi_{i, \vec{\alpha}} = \vec{x}_i^{\vec{\alpha}}$ is the \textit{polynomial matrix}, and
$(\vecf_E)_i = f(\vecx_i)$ for $\vecx_i \in S(E)$.

Before turning to error estimation and convergence theory, we comment on key
parameters in the RBF interpolant.
We recommend first choosing $\textrm{deg}$, the degree of the
appended polynomial, which governs the order of accuracy as $h \to 0$. 
This choice determines the size of the polynomial matrix $\Pi$, which must have
more rows than columns to ensure the system \cref{eqn:interp_system} is invertible. 
Specifically, the stencil size must satisfy $k \ge \genfrac(){0pt}{1}{\textrm{deg}+\text{dim}}{\textrm{deg}}$. 

We will specify our choice of $k$ for all numerical experiments, though there are often reasons to choose $k$ larger than the minimum.
As general principle, increasing $\textrm{deg}$ improves accuracy, while increasing $k$
enhances stability. For instance, near domain boundaries, stencils may become
one-sided, potentially reducing accuracy due to Runge's
phenomenon. Enlarging stencils in such regions can mitigate this effect and improve accuracy without changing $h$. For further discussion on the interplay between stencil size and polynomial degree, see~\cite{FLYER2016, BAYONA2017, BAYONA2019}.


\rev{
We must also choose $\ell$, the degree of the polyharmonic spline (see
\cref{eq:polyHarmonic}). Although increasing $\ell$ can reduce error, the 
improvement diminishes quickly, and non-singularity of the interpolation 
matrix requires $\text{deg} \ge \lfloor \ell /2 \rfloor + 1$ 
\cite{WENDLAND2005, FASSHAUER2007}. In practice, we therefore adopt 
parameters that balance accuracy, stability, and cost: stencil sizes 
$k=21$ or $32$, slightly larger than the minimum needed for unisolvency, 
improve stability near boundaries without substantially increasing cost 
\cite{FLYER2016, BAYONA2017}; and we fix $\ell=3$, since larger values 
yield little accuracy gain but worsen conditioning \cite{WENDLAND2005}. 
With appended polynomials up to degree $\deg=4$, these choices consistently 
performed well, and convergence tests (Sections~\ref{ssec:rbfqf},~\ref{ssec:manifold}) 
confirm robust accuracy across both flat and curved geometries. 
}

We conclude this section by interpreting the interpolant $s$ defined in
\cref{eq:sFunDef} as a projection $P_n f$, aligning this step with the framework in \Cref{sec:model}. While the projection
operator $P_n: C(D) \to C(D)$ is expected to preserve continuity, the interpolant $s$ is constructed locally and may be discontinuous across element
boundaries. A continuous projector $P_n$ can be obtained by applying weighted
averages to neighboring local interpolants.

For example, consider a Voronoi partition $\{V_i\}_{i=1}^n$ of $D$, and
let $\{T_{i, j, k}\}$ denote the associated Delaunay triangulation, where $1 \le i,
j, k \le n$ index the vertices of the triangle. 
For any $\vecx \in T_{i,
j, k} \subseteq D$, let $b_{i, \vecx}, b_{j, \vecx}, b_{k, \vecx}$ be the barycentric
coordinates of $\vecx$ with respect to $\vecx_i, \vecx_j, \vecx_k \in \Xi$. We define the piecewise projector
\[
    (P_n f)(\vecx) = \sum_{T_{i, j, k}} \chi_{T_{i, j, k}}(\vecx) \bigg(
        b_{i, \vecx} s_{V_i}(\vecx) + 
        b_{j, \vecx} s_{V_j}(\vecx) + 
        b_{k, \vecx} s_{V_k}(\vecx)
    \bigg),
\]
which is continuous across $D$ and smooth within each triangle. With $P_n$ defined, the Lagrange basis $\{L_1, \dots, L_n\}$ is given by $L_i = P_n v_i$ where each $v_i \in C(D)$ satisfies $v_i(\vecx_j) = \delta_{ij}$.

Second, the projection method requires $P_n f \to f$ as $n \to \infty$ but the
relationship between $n$ and our interpolant is not apparent. To relate them, we
consider a sequence of interpolation node sets $\{\Xi_n\}$ with covering radii
$h(\Xi_n) \to 0$ as in \cref{eqn:spacing}. Such sequences can be easily constructed by refining triangulations or by placing points to efficiently fill
the domain.

This framing aligns with interpolation convergence theorems in the literature
which establish error bounds for polyharmonic spline interpolation such as
\[
    \|f - P_n f\|_\infty \le C h^{\textrm{deg}}
\]
or better in some cases~\cite{WENDLAND2005, FASSHAUER2007, buhmann2000, johnson1998}.
These results apply to \textit{global} interpolation ($k=n$). 
While we found no published proof that \textit{local}
interpolation achieves the same rate, we expect such a result could follow from adaptations of the existing global proofs~\cite{FASSHAUER2007,
WENDLAND2005}, though it would likely require geometric justification that the cone
condition can be relaxed, except near boundaries. 
In practice, studies using local interpolation offer numerical evidence of convergence, often showing convergence rates faster than $O(h^\textrm{deg})$~\cite{FLYER2016, BAYONA2017, BAYONA2019} which we observe in \Cref{sec:numerical_tests}.

\subsection{Quadrature In Flat Domains} \label{ssec:rbfqf}
The theory of RBF-QF has its roots in finite difference methods,
known as RBF-FD. 
Both quadrature and differentiation are linear functionals -- linear maps into $\mathbb{R}$ -- and share a common theoretical foundation. 
Derivative approximation using RBFs was introduced first~\cite{KANSA1990}, and the modern
RBF-FD method based on local interpolation was discovered independently by
several groups shortly after~\cite{wright2003radial, tolstykh2003using, shu2003local, cecil2004numerical}. 
Due to this history, some authors refer to RBF-QF as a variant of RBF-FD, though it produces quadrature weights rather than finite difference weights. 
We begin with the general theory for approximating linear functionals,
then specialize to quadrature. 
This approach clarifies the distinction between RBF-QF and
RBF-FD, and highlights when and why a mesh is needed.

To approximate a linear functional $\LL$ (such as quadrature or differentiation at a
point) applied to a function $f: D \to \RR$, we first consider the local RBF
interpolant $s \colon D \to \RR$ which, by \cref{eq:sFunDef}, gives
\begin{align*}
    \LL f &\approx \LL s 
        = \sum_{E} \LL  s_E \chi_{E}
        = \sum_{E} \sum_{\vecy \in S(E)} c_{E,\vecy} \LL ( \phi_{\vecy}\chi_E \big) +
        \sum_{|\vec{\alpha}| \le \text{deg}} d_{E, \vec{\alpha}} \LL[ \chi_E
        \pi_{\vec{\alpha}} ].
\end{align*}
The coefficients $c_{E, \vecy}$ and $d_{E, \vec{\alpha}}$ are given in vector form in \cref{eqn:interp_system}. 
Representing the previous equation as a sum over dot products, we have
\[
\begin{aligned}
    \LL f &\approx \sum_E \begin{bmatrix} \LL \vec{\phi} 
      & \LL \vec{\pi} \end{bmatrix} \begin{bmatrix} \vec{c} \\ \vec{d} \end{bmatrix} 
      &&= \sum_E \begin{bmatrix} \LL \vec{\phi} & \LL \vec{\pi} \end{bmatrix}
        \begin{bmatrix} A & \Pi \\
        \Pi^T & 0
        \end{bmatrix}^{-1}
        \begin{bmatrix} \vec{f}_E \\ \vec{0} \end{bmatrix}, \\
    &= \sum_E 
        \begin{bmatrix} \vec{w}_E & \vec{\gamma}_E \end{bmatrix} 
        \begin{bmatrix} \vec{f}_E \\ \vec{0} \end{bmatrix} 
    &&= \sum_{\vecx_i \in \Xi} w_{\vecx_i} f(\vecx_i), \\
    &= \vecw^T \vecf, 
    &&
\end{aligned}
\]
where $\vec{w}_E$ denotes the element-wise functional weights, and $\vec{w}$ the combined functional weights,
\begin{equation} \label{eqn:functional_stencil}
    \vec{w}_E =
        \begin{bmatrix} A & \Pi \\
        \Pi^T & 0
        \end{bmatrix}^{-1}
        \begin{bmatrix} \LL\vec{\phi} \\ \LL\vec{\pi} \end{bmatrix},
        \qquad 
    (\vecw)_{\vecx_i} = w_{\vecx_i} = \sum_{E \colon \vecx_i \in S(E)} (\vec{w}_E)_{\vecx_i}. 
\end{equation}
If $\LL$ is a differential (integral) operator, then $\vec{w}$ represents finite
difference (quadrature) weights. Moreover, the expressions involving $w_{\vec{x}_i}$ and
$\vec{w}$ help outline the implementation of the algorithm.

First, for each element $E$, we identify the stencil $S(E)$, and evaluate the functionals
$\LL (\phi_{\vecx_i}\chi_E) $ and $\LL (\vecx^{\vec{\alpha}} \chi_E)$.
We then form the system in \cref{eqn:functional_stencil} and solve for the
weights $\vecw_E$, whose order reflects the ordering of the nodes in the stencil
$S(E)$. Each node $\vecx_i$ may appear in multiple stencils and thus be assigned
several weights $(\vecw_E)_{\vecx_i}$. The final quadrature weights are obtained by
summing all contributions associated with each node across stencils.

Since we focus on quadrature, we set 
$\LL f = Qf = \int_D f(y) d\mu(y)$ as in \Cref{sec:model}. While our theory suggests
$Q_n = QP_n$, we use the simpler form $Q_n = Qs = \vec{w}^T \vec{f}$, which offers
the same accuracy and is easier to compute. Each functional evaluation $\LL
\vec{\varphi}, \LL\vec{\pi}$ thus involves analytically integrating an RBF or polynomial
basis function over an element $E$. The resulting quadrature approximation is then
$Q_n f = \vec{w}^T \vec{f}$.

\subsection{Quadrature on Smooth Closed Manifolds} \label{ssec:manifold}
When the cortical domain $D$ is a $2$-manifold with curvature, applying RBF-QF
is more involved. Convergence theory for such schemes remains an active area of research~\cite{glaubitz2023towards}, while implementations have been developed by Reeger and Fornberg~\cite{REEGER2016Sphere, REEGER2016,
REEGER2018}. We summarize their approach for smooth closed surfaces and refer
to~\cite{REEGER2018} for the treatment of surface boundaries. Numerical evidence of convergence rates is presented in \Cref{ssec:surface_convergence}.

\begin{figure}
    \centering
    \includegraphics[width=\textwidth]{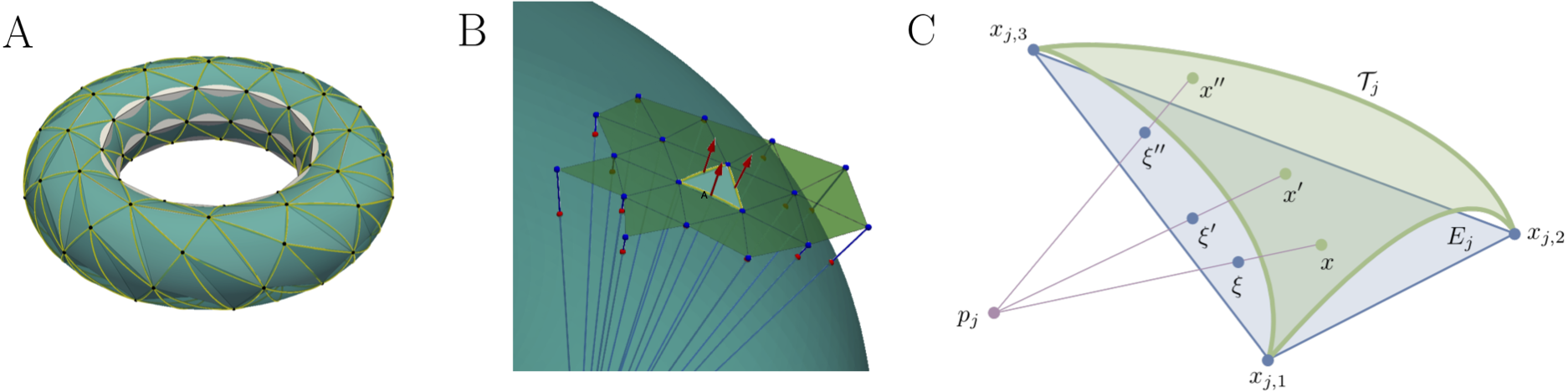}
    \caption{
        \textbf{A.} The surface of a torus ($D_T$, transparent light-green) is sampled at
        some well-distributed points ($\Xi$, black) which are then used for a triangulation
        ($\{E_j\}_{j=1}^M$, white flat triangles). 
        Using the method described in
        \Cref{ssec:manifold}, these planar triangles are projected onto the torus to form a
        partition of surface triangles ($\{\mathcal{T}_j\}_{j=1}^M$, outlined in yellow).
        \textbf{B.}~Stencil and projected stencil of a flat triangle ($E_j$, white, but appears light green because it is behind the surface) from a surface
        triangulation. Nearby surface points form the stencil ($S(E_j)$, red) and are projected along lines (blue) emanating
        from the projection point, onto co-planar points ($\xi$, blue). The projection point is computed from the vertices of the
        triangle and approximated edge normals (red arrows).
        \textbf{C.} A detailed diagram relating a
        flat triangle $E_j$ and its associated surface triangle $\mathcal{T}_j$.
        Points $\vec{x} \in \pi_j$ on the plane of $E_j$ are
        projected onto $\vec{\xi} \in \mathcal{T}_j$ along the lines
       through the projection point $p_j$. 
    }
    \label{fig:torus_partition}
    \label{fig:sketch-tri}
\end{figure}

The method begins with a set of interpolation nodes $\Xi$ and defines two sets of elements. 
The first is a planar triangulation $\calE = \cup_{j=1}^m E_j$ with nodes exactly equal to $\Xi$. 
The method is designed for curved (non-planar) domains where the cortical surface $D$ does not coincide with $\calE$. It constructs a second set of elements forming a partition $D = \cup_{j=1}^m \calT_j$, where each
$\calT_j$ is a curved triangular patch on $D$, in contrast to the flat triangles $E_j$ in $\calE$ (See \cref{fig:sketch-tri}). The two families
share the interpolation nodes, $\Xi \subseteq D \cap \mathcal{E}$, and are paired such that for each $j$ there exists a bijection $\mu_j(E_j) = \calT_j$.
We defer a precise definition of the mappings to a later section, and simply note
that once the families $\{ E_j \}$, $\{ \calT_j \}$, and $\{ \mu_j \}$ are in place, one
can write, without approximation,
\begin{equation}\label{eq:integral}
  \int_{D}f(\vecx) \,dS = \sum_{j=1}^{m} \int_{\calT_j} f(\vecx) \,dS 
  = \sum_{j=1}^{m} \int_{E_j} f\bigl(\mu_j(\vec{\xi})\bigr) \, 
  \| (\partial_{\xi_1} \mu_j \times \partial_{\xi_2} \mu_j)(\vec{\xi}) \|\,d\vec{\xi},
\end{equation}
in which 
$
  \| (\partial_{\xi_1} \mu_j \times \partial_{\xi_2} \mu_j)(\vec{\xi}) \|
$
is the Jacobian of $\mu_j$ at $\vec{\xi}$. The method introduced by Reeger and
Fornberg amounts to~\cite{REEGER2016,REEGER2016Sphere,REEGER2018}: (i) defining the mappings $\{ \mu_j \}$; (ii) selecting a quadrature scheme for the integral over $E_i$, on the right-hand side of \cref{eq:integral}.

The mappings $\{  \mu_j \}$ depend on a set of \textit{projection points} $\{
\vec{p}_j \}$, as shown in \Cref{fig:sketch-tri}. For now, we assume these
points are given and explain later how to determine them. For each fixed
$j$, let $\pi_j$ denote the plane in which the triangle $E_j$. We define
\[
  \mu_j \colon \pi_j \to D, \qquad \vec{\xi} \mapsto \vec{x}
\]
where $\vec{x} \in D$ is the closest point to $\vec{\xi} \in \pi_j$ along the line connecting
$\vec{p}_j$ and $\vec{\xi}$. Such a point exists provided $\vec{\xi}$ $\vecx$ lies in a sufficiently small neighborhood of $E_j$ and the node set is $\Xi$ is dense enough. The mapping $\mu_j$ then
sends $E_j$ to a curved triangle $\calT_j \subseteq
D$. This construction explains why $\{\vec{p}_j\}$ are called projection points: if
$\vec{p}_j$ were a light source, then $\calT_j$ would be the shadow of $E_j$ cast onto the
curved surface $D$. The point $\vec{x} \in \calT_j$ is defined as the closest such intersection since the projection line may
intersect $D$ more than once. Finally, we note that $\mu_j$ is defined for
points on $\pi_j$ even outside $E_j$.

By construction, $E_j$ and $\calT_j$ intersect at the nodes $\{ \vecx_{j,1}, \vecx_{j,2},
\vecx_{j,3} \} \subset \Xi$ (See \cref{fig:sketch-tri}). The mapping $\mu_j$ sends each edge of the
planar triangle $E_j$ connecting $\vecx_{j,k}$ and $\vecx_{j,l}$ to the corresponding edge of the curved patch $\calT_j$, preserving endpoints. Each such edge, together with
the projection point $\vec{p}_j$, lies in a common \textit{cutting plane}; thus
$\vec{p}_j$ is the intersection of three cutting planes.

The definition of the mappings $\{ \mu_j \}$ and the properties above are valid
for any choice of projection points $\{ \vec{p}_j \}$. However, these points must be selected carefully to ensure that $\{\mathcal{T}_j\}$
forms a valid partition of $D$. 
Reeger and Fornberg outline a method to choose projection points based on the
normal vectors $\{\vec{n}_j\}$ of the planar triangules $\{E_j\}$. 
For any pair of adjacent triangles $E_j$ and $E_{j'}$, sharing an edge $E_j \cap
E_{j'}$, a necessary condition for $\{ \calT_j \}$ to partition $D$ is
that the maps $\mu_j$ and $\mu_{j'}$ both send the shared edge to the same
curved edge in their respective curved triangles, that is
\[
  \mu_j(E_j \cap E_{j'}) = 
  \calT_j \cap \calT_{j'} = \mu_{j'}(E_j \cap E_{j'}),
  \qquad 
  \text{for all $j$ and $j'$,}
\]
which implies that $\vec{p}_i$ and $\vec{p}_j$ must lie in the same cutting plane.
They resolve this by choosing the cutting plane orthogonal to $\vec{n}_j - \vec{n}_{j'}$ with the condition $\vec{p}_j \cdot (\vec{n}_j - \vec{n}_{j'}) = 0$. This approach determines the distribution $\{ \vec{p}_j \}$, as each planar triangle is associated with three such cutting planes, uniquely fixing its projection point.

Once the bijections $\{ \mu_j \}$ are defined, a quadrature rule for integrands $g
\colon E_j \to \RSet$ can be written as
\[
\int_{E_j} g(\vec{\xi}) \,d\vec{\xi} \approx \sum_{r=1}^{q_j} g(\vec{\xi}_{j,r})
  \rho_{j,r},
\]
where the quadrature nodes satisfy
\[
  \mu_j(\vec{\xi}_{j,r}) = \vec{x}_{j,r}, \qquad r \in \NSet_{q_j}, \qquad j \in \NSet_m,
\]
and $ \{ \vec{x}_{j,1}, \ldots, \vec{x}_{j,q_j} \} = S(E_j)$; that is, the quadrature nodes
for $E_j$ are the pre-images of the stencil points under $\mu_j$ (see \Cref{ssec:local_interpolation} for a definition of $S(E_j)$). Importantly, the nodes $\{\vec{\xi}_{j,r} \}_{r =1}^{q_j}$ lie in the plane $\pi_j$, but not necessarily within the triangle $E_j$ itself, as illustrated in \Cref{fig:sketch-tri}.
The points $\mu_j^{-1}(S(E_j)) \in \pi_j$ now form a planar stencil for $E_j$, and we use the procedure outlined in \Cref{ssec:rbfqf}. Specifically, we treat each \(E_j\) as a flat domain and apply RBF-QF locally, using the basis functions evaluated at the stencil nodes and integrating them over \(E_j\). This yields weights that, when paired with function evaluations on the surface at \(\vec{x}_{j,r} = \mu_j(\vec{\xi}_{j,r})\), define a quadrature rule that approximates integrals over the curved domain \(D\).

\section{Numerical Experiments for Neural Fields} \label{sec:numerical_tests}

This section presents numerical experiments using RBF-QF, both as a standalone quadrature method and within a method of lines neural field simulation. 
We begin by introducing notation we will use throughout the section.

We consider three spatial domains of integration: (i) $D_1 = [0, 1]^2$, the
unit square; (ii)
$D_{2\pi} = [-\pi, \pi]^2$ the square of width $2\pi$ centered at the
origin (iii) $D_T$, a torus with major radius $R = 3$ and minor radius of $r=1$,
parameterized by angles $(\phi, \theta) \in D_{2\pi}$
as
\begin{align}
    \vecx(\phi, \theta) &= \begin{bmatrix}
        x(\phi, \theta) \\
        y(\phi, \theta) \\ 
        z(\phi, \theta)
    \end{bmatrix}
    =
    \begin{bmatrix}
        (R + r\cos \theta)\cos \phi  \\
        (R + r\cos \theta)\sin \phi  \\
        r\sin \theta
    \end{bmatrix}. \label{eqn:torus}
\end{align}
We denote the inverse map by $\psi: D_T \to D_{2\pi}$.

Unless otherwise noted, we use the standard Euclidean distance metric in $\RR^2$ or
$\RR^3$, as appropriate, for stencil selection and evaluating basic function inputs.
For certain test functions and neural field solutions on $D_1$ or $D_{2\pi}$, we
instead use a doubly periodic distance metric with periods $\Lambda = 1$ and $
\Lambda = 2\pi$,
respectively, defined by
\begin{align}
    \text{dist}(\vecx, \vecy) = \min_{i,j=-1, 0, 1} \left\| \vecx - \vecy + \Lambda \left( i\begin{bmatrix}1 \\ 0 \end{bmatrix} + j\begin{bmatrix} 0 \\ 1 \end{bmatrix} \right) \right\|. \label{eqn:dist}
\end{align}

Most of our test functions, as well as both neural field solutions, use scaled Gaussian functions defined with respect to this periodic distance:
\begin{align}
    \text{Gauss}(\vecx, \vecy; \sigma)
        &= \frac{1}{2\pi \sigma^2}\exp\left[- \frac{1}{2\sigma^2} \ \text{dist}^2(\vecx, \vecy)\right]. \label{eqn:gauss}
\end{align}
In this setting, the kernel function is defined as $w(\vecx, \vecy) =
\tilde{w}(\text{dist}(\vecx, \vecy))$, where $\tilde{w}$ is a standard radial
function (e.g., Gaussian). While $\tilde{w}$ is not periodic, periodicity is enforced
via the modified distance function in \cref{eqn:dist}, which wraps Euclidean distance
across boundaries—effectively defining $w$ on a flat torus. This approach ensures
that $w$ depends on the shortest periodic path between $\vecx$ and $\vecy$, not just
their Euclidean separation.

In neural field simulations we use the sigmoidal firing rate function
\begin{align}
    f[u] &= \bigg(1 + \exp\big(-\gamma(u - \vartheta)\big)\bigg)^{-1} \label{eqn:firing_rate}
\end{align}
with gain $\gamma = 5$ and threshold $\vartheta = 1/2$.

To test the quadrature within a method of lines neural field simulation, we proceed
with a {\em method of manufactured solutions} to construct equations with known analytic
solutions. Specifically, we choose a solution $u(t, \vecx)$ and define a forced
neural field equation parameterized by a weight kernel $w$ and the firing firing-rate
function $f$: 
\begin{align}
    \partial_t u(t,\vecx) = -u(t,\vecx) + \int_D w(\vecx, \vecy) f[u(\vecy)] \ dy + F(t, \vecx) \label{eqn:manufactured_system}
\end{align}
where the forcing term $F$ is chosen to ensure that $u$ satisfies the equation exactly.

\subsection{Node selection}
We briefly describe our selection of quadrature nodes and triangular meshes in each
domain. Rather than detailing each construction, we provide simplified description
and refer interested readers to the Code availability section (\cref{sec:code})
for implementation details, including quadrature weight generation, node
placement, and time integration.
For the square domains $D_1$ and $D_{2\pi}$ we compare two node families. The first
places interior nodes at the vertices of a regular equilateral triangle tiling, with
equally spaced nodes along the boundary. We refer to this as the {\em regular grid},
since interior stencils are rotated translates of each other and yield identical
weights (\Cref{fig:flat_random_quad_error}E).
The second family distributes random interior nodes using a deterministic
point-repulsion process, again placing equally spaced nodes along the boundary. An
example is shown in \Cref{fig:flat_random_quad_error}C.
In both cases, we use Delaunay triangulation of the quadrature nodes to define the
mesh.

For quadrature nodes on the torus $D_T$, we use a deterministic regular grid of {\em
spiral nodes}. They are defined by a triangular tiling of the parameter space $D_{2\pi}$ that avoids thin triangles when mapped onto the torus.
The nodes lie on lines of constant $\theta$ (rotation
around the minor axis) and are equally spaced in the $\phi$-direction, but staggered
so they form approximately $60^\circ$ angles with the $\theta$-constant lines. These angled lines in trace closed spiral paths on the torus surface. While the
tiling is regular in parameter space, it does not account for the curvature of the surface, and thus triangle sizes vary—becoming larger in regions of positive curvature and smaller in regions of negative curvature.

\subsection{Quadrature Experiments on the Unit Square}\label{ssec:square}

It is standard to verify quadrature convergence by measuring relative error as the
mesh is refined (i.e., as the spacing $h \to 0$). We have performed numerous such
experiments and present a representative sample in \Cref{ssec:flat_convergence}.
Before turning to convergence rates, we first examine specific meshes to illustrate
how their structure influences the resulting quadrature weights and errors. This
focus is worthwhile because unlike conventional quadrature methods, which require
specific node placements (e.g., Gaussian or Newton-Cotes), RBF-based quadrature can
accommodate arbitrary node sets. In this setting, mesh spacing $h$ is a useful
summary statistics but does not fully determine the node layout or resulting weights.

\begin{figure}
    \centering
    \includegraphics[width=0.95\textwidth]{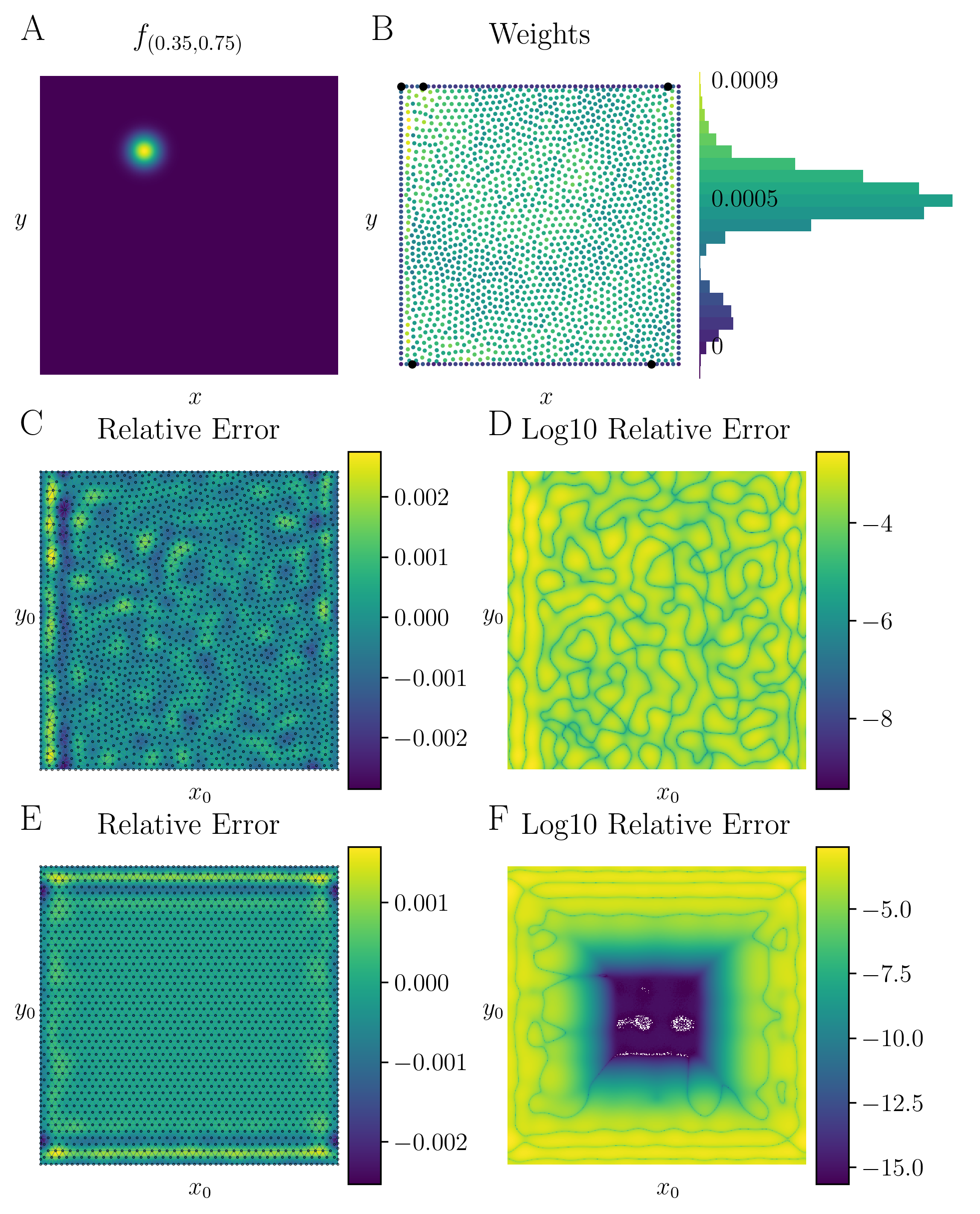}
    \vspace{-3mm}
    \caption{\textbf{RBF-QF on the Unit Square.}
        \textbf{A.} A single test function as in
        \cref{eqn:periodic_gaussian_test_function}).
        \textbf{B.} Quadrature node locations for $n = 2000$ points, colored by
        their associated weights. Boundary-adjacent nodes have smaller weights due to
        asymmetric stencils, and five nodes yield negative weights (black circles).
        Right: a histogram of weights showing a bimodal distribution—one mode near
        $1/n = 0.0005$ and a second associated with boundary clustering.
        \textbf{C–D.} Relative error (\textbf{C}) and log$_{10}$ relative error
        (\textbf{D}) of the quadrature rule applied to Gaussian test functions
        $f_{\vec{x}_0}$ centered throughout the domain. 
        \textbf{E–F.} Same as panels \textbf{C–D}, but using a near-regular
      triangular grid. Interior weights approach $1/n$ and the error approaches
    machine precision away from the boundary.}
    \label{fig:flat_random_quad_error}
    \vspace{-8mm}
\end{figure}

To demonstrate more precisely, we introduce a qualitative test in which a localized
test function smooths the quadrature operator.
Specifically, we define
a family of steep Gaussian test functions $f_{\vecx_0}: D_1 \to \RR$ by
\footnote{
  The prefactor of $2 \pi 10^2$ has been used in the numerical experiment, but it has
no impact on the relative error presented below, and can be safely ignored}
\begin{align}
    f_{\vecx_0}(\vecx) & = 2 \pi 10^2 \; \text{Gauss}(\vecx, \vecx_0; 10) \label{eqn:periodic_gaussian_test_function}
\end{align}
where each function is centered at $\vecx_0 \in D_1$. 
These Gaussians, defined in \cref{eqn:gauss}, using the doubly-periodic distance
\cref{eqn:dist}, integrate to the same constant regardless of their center. The mass
is sharply localized near $\vecx_0$, with rapid decay in all derivatives away from
the center. While the periodic distance is not smooth, the resulting test functions
are effectively smooth to machine precision. An example is shown in
\Cref{fig:flat_random_quad_error}A.

\Cref{fig:flat_random_quad_error}B shows the quadrature nodes, colored by their
associated weights. We observe that weights along the boundary are consistently
smaller. This is partly because boundary nodes appear in fewer stencils, and partly
due to mild clustering near the edges—an intentional design to mitigate
boundary-related errors such as the Runge phenomenon. Five nodes have negative
weights (highlighted with black circles), though each is small in magnitude. The
histogram to the right reveals a bimodal distribution: the larger mode is centered
near $1/n = 0.0005$, consistent with the average area per node, while the smaller
mode corresponds to the boundary-adjacent nodes. This empirical structure reflects
both the local stencil geometry and the rule’s enforcement of exactness for constant
functions.

We generate approximately $n=2000$ random and regular nodes (as
described in \Cref{sec:numerical_tests}) and apply the RBF-QF algorithm from
\Cref{ssec:rbfqf} to compute quadrature nodes $\Xi$ and associated weights
$\{w_{\vec{\xi}}\}_{\vec{\xi} \in \Xi}$.
To assess spatial error, we introduce a relative error as a function of the test
function center:
\begin{equation}
  E(\vecx_0) = 
  \frac{Qf_{\vecx_0} - Q_n f_{\vecx_0}}{Q f_{\vecx_0}}
  \label{eqn:flat_spatial_error}
\end{equation}
Large values of $|E(\vecx_0)|$ indicate substantial quadrature error near $\vecx_0$,
suggesting the local node configuration contributes disproportionately to the total
error. This can be interpreted as measuring the relative error in approximating the
convolution $(1 \ast f_{(\cdot)})(\vecx)$
via our quadrature rule.

We use localized Gaussian test functions to visualize how quadrature error varies
spatially for different node configurations. In \Cref{fig:flat_random_quad_error}, we
demonstrate this on two node sets in the unit square
(See Section \ref{sec:numerical_tests}).
\Cref{fig:flat_random_quad_error}\textbf{C} shows a heatmap of $E(\vecx_0)$,
the relative signed quadrature error of our Gaussian test functions $f_{\vecx_0}$. We
again use $n=2000$ quadrature nodes, the basic function $r^3$, a stencil size of
$k=21$, and append third-degree polynomials. The error is relatively small, varies
continuously, and is not systematically positive or negative. Most of the error
magnitudes are well below the extremal values, highlighting the importance of testing
across a variety of test functions.

For this particular mesh, the points are roughly evenly spaced except near
the top and bottom boundary where they cluster. Stencils near the boundary are
necessarily one sided which can lead to Runge phenomenon \cite{FLYER2016, fornberg2007runge}.
Clustering nodes near these boundaries is one way to combat this and demonstrates the
utility of our spatially localized test functions. We observe that the error is
generally higher along the left and right boundaries, where no clustering is applied,
and comparatively lower near the top and bottom edges, where node clustering helps
control boundary effects.

Any test function can exhibit surprisingly low error due to the random nature
of the quadrature nodes. The function $E(\vecx_0)$ is continuous and oscillatory, and
since the quadrature rule integrates constants exactly, the average of $E$ over the
domain is zero. The continuity of $E(\vecx_0)$ follows from the fact that $Q_n f_{\vecx_0}$ is a continuous function of $\vecx_0$, and $Qf_{\vecx_0}$ is constant in $\vecx_0$. The average value of $E$ is zero because the quadrature weights are exact for constant functions, a fact that can be verified by a short computation. 
Thus,
there must exist closed curves along which $E(\vecx_0) = 0$. We
visualize this zero-level set by plotting $\log|E(\vecx_0)|$ in
\Cref{fig:flat_random_quad_error}D. While the overall structure remains similar
across different meshes, the precise locations of high and low error regions will
vary. This implies that, even for a fixed test function, quadrature error exhibits
mesh-dependent spatial variation—even when meshes have similar $n$ and point
densities. For this reason, it is common in RBF quadrature studies to generate
multiple random node sets for each choice of mesh parameters.

In practice, the presence of small negative weights does not prevent convergence. As
noted earlier, \Cref{fig:flat_random_quad_error}B shows five such weights, each small
in magnitude, and their influence appears minimal in
\Cref{fig:flat_random_quad_error}C–D. It is often proposed that a quadrature rule is
stable provided
\[
    \sum_{n} w_{\vec{\xi}_n} - \sum_{n} |w_{\vec{\xi_n}}| = 0
\]
i.e., if all weights are non-negative \cite{glaubitz2023towards}. High-order Newton–Cotes rules, for instance,
fail this criterion and are indeed unstable. However, this notion of
stability seems overly restrictive for our purposes. A more lenient condition--—one that
permits convergence—--is that the difference between the total signed and absolute
weights remains bounded as $n \to \infty$ \cite{trefethen2022exactness,
polya1933konvergenz}. Indeed, the RBF-QF method has been analyzed under this
framework \cite{glaubitz2023towards}, and our experiments confirm that it converges
robustly even in the presence of small negative weights. 
\rev{
Alternative strategies,
such as $L^1$ optimization of stencil weights~\cite{abrahamsen2019explicit}, have been
proposed to mitigate or eliminate negative weights while maintaining accuracy, but we
do not pursue these modifications here.
}

Lastly, Panels E and F of \Cref{fig:flat_random_quad_error} illustrate a case where
the quadrature rule achieves spectral accuracy away from the boundary. The test functions and algorithmic
parameters remain unchanged, but the nodes are now arranged on a near-regular
triangular grid. In this symmetric configuration, the error is extremely small away
from the boundary and reaches machine precision near the domain center. The white
pixels visible in Panel F result from floating-point coincidence: the analytic
integral and quadrature evaluation yield exactly the same double-precision value.

This level of accuracy arises from geometric regularity. When the mesh elements and
their associated stencils are rotated translates of each other, the resulting
quadrature weights are identical. Coupled with exact integration of constant
functions, this symmetry forces all interior weights to equal the average area per
node. Importantly, this value is independent of the degree of appended polynomial
basis terms, which determine the formal convergence rate. As a result, the quadrature
rule achieves arbitrarily high order in the interior—analogous to the spectral
accuracy of the trapezoidal rule for periodic functions in one dimension. \rev{We observe the same phenomenon on other geometries: when nodes are placed on
regular, nearly uniform grids on the square (Fig.~SM1), sphere (Fig.~SM2), or cyclide (Fig.~SM3), the quadrature
achieves high accuracy, as documented in the Supplementary Materials.}

\subsection{Convergence on the unit square}\label{ssec:flat_convergence}

\begin{figure}
    \centering
    \includegraphics[width=1.0\textwidth]{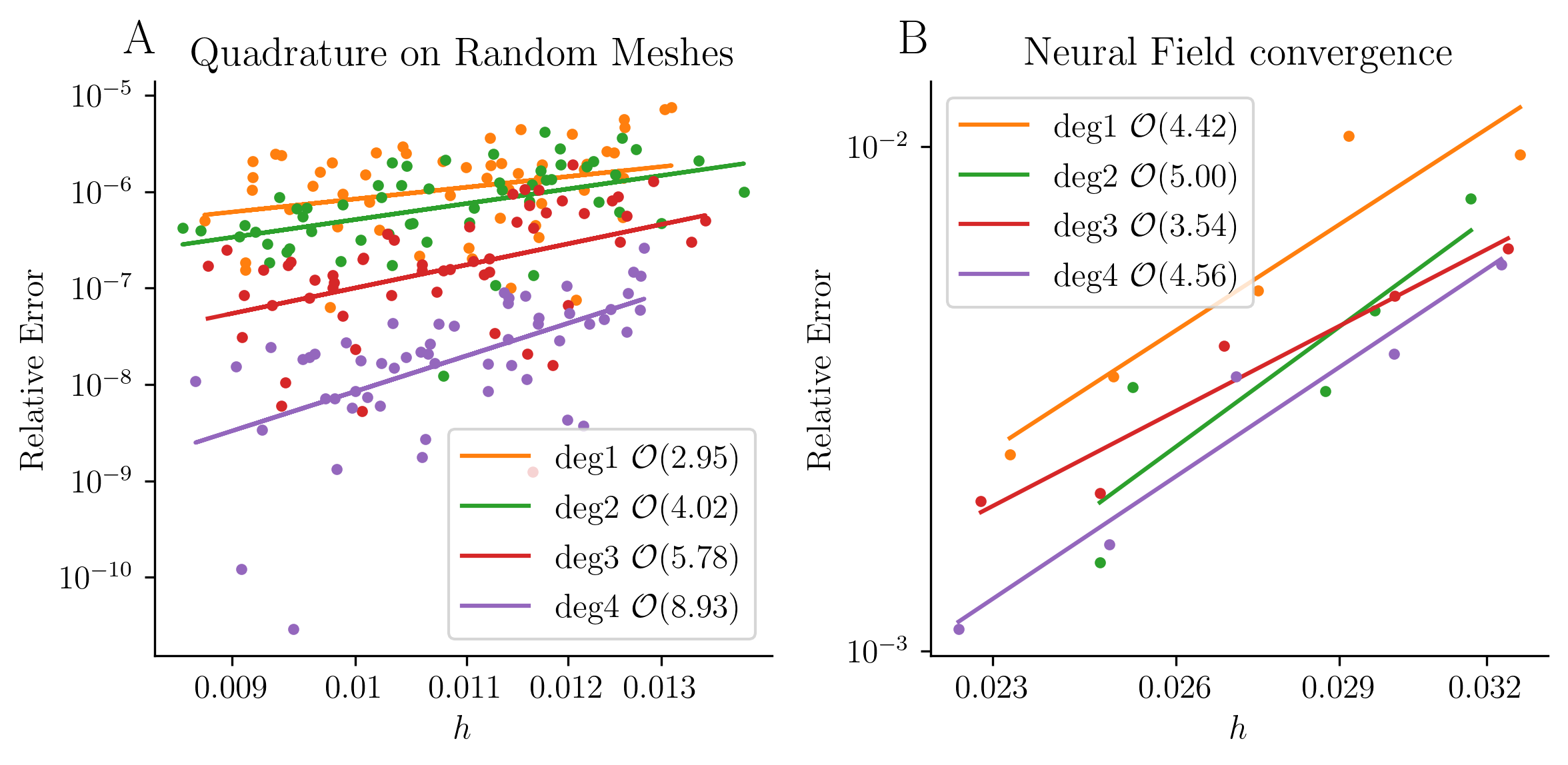}
    \caption{
        \textbf{A}~Convergence of RBF-QF quadrature for the function
        $T_5(2x-1)T_4(2y-1)+1$ on the unit square. 
        We use the RBF $\phi(r) = r^3$, stencil size $k=21$, and
        append polynomials terms up to the specified degree (denoted by color).
        \textbf{B}~Convergence of a neural field simulation using the same quadrature rules (see main text for details).
    }
    \label{fig:flat_convergence}
\end{figure}

We evaluate convergence of RBF-QF and neural field simulations
on the unit square, using randomly chosen quadrature nodes and triangulations (See
\Cref{sec:numerical_tests}). The test function is $f(x, y) = T_5(2x - 1) T_4(2y-1)$,
a product of Chebyshev polynomials (\Cref{fig:flat_convergence}~C). Using
$\phi(r) = r^3$ and stencil size $k=21$, we compute quadrature weights. As shown in
\Cref{fig:flat_convergence}A, errors decrease rapidly with increasing
appended polynomial degree, often faster than $\mathcal{O}(h^\text{deg})$.
Panel~{B} shows similar convergence for the simulation of the neural field
solution, which we now describe.

We next present measured convergence rates for a method of lines simulation of a neural field, using RBF-QF to discretize the weight kernel. 
The manufactured solution is
\begin{align*}
    u(t, \vecx) &= 
        f^{-1}\left[ \text{Gauss}(\vecx, \vecx_0(t); \sigma=11/10) +
        1/10\right], \\
    \vecx_0(t) &= [\cos(t)/5, \sin(t)/5]^T,
\end{align*}
which is indeed a solution to \cref{eqn:manufactured_system} on $(t, \vecx) \in
[0, 1/10] \times D_{2\pi}$.
Here $f$ is defined in \cref{eqn:firing_rate}, and the weight kernel is $w(\vecx,
\vecy) = \text{Gauss}(\vecx, \vecy; \sigma=1/40)$ (See \cref{eqn:gauss}).
These expressions are used to construct a consistent forcing function.
\Cref{fig:flat_convergence}
{B} reports convergence results for different polynomial degrees. In each case, the
observed order of accuracy exceeds the expected $\mathcal{O}(h^\text{deg})$ rate. 
\rev{
For random meshes, the noise in the error curves is expected: each random placement
of nodes slightly alters local point density and interpolation accuracy, and these
local variations accumulate in the global error measure, producing the observed
fluctuations.
}

\begin{figure}
    \centering
    \includegraphics[width=0.95\textwidth]{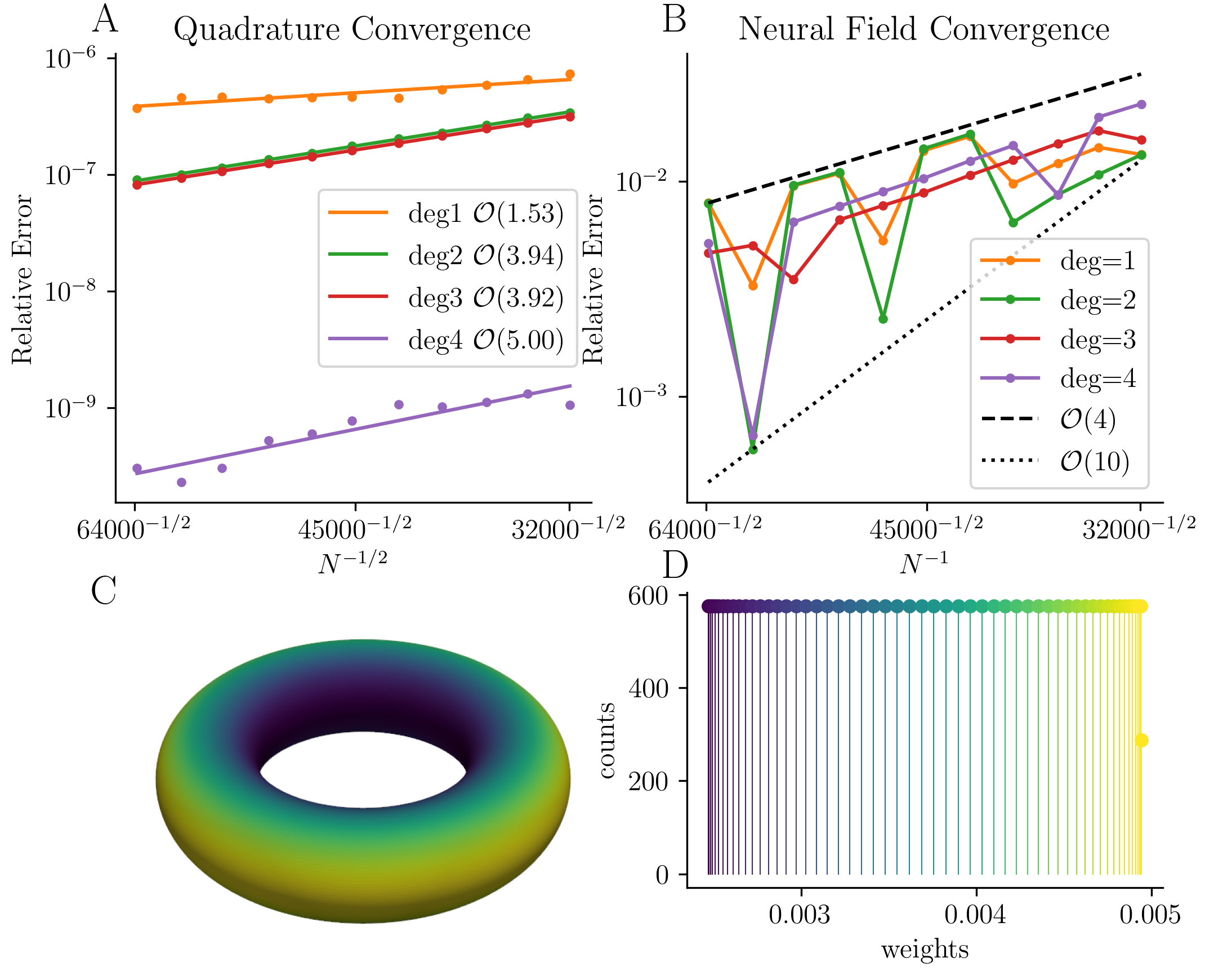}
    \caption{
        \textbf{A.}~Convergence of surface quadrature on the torus $D_T$ using the test function $f(x, y, z) = \sin(7x)+1$. In all cases, observed rates exceed the degree of appended polynomials.
        \textbf{B.}~Convergence of a neural field simulation using the same surface quadrature rule.
        \textbf{C.}~Torus surface colored by the quadrature weights for $k=12$ and polynomial degree 2. Due to mesh regularity, weights are constant along $\theta$ and larger in regions of positive curvature where triangles are bigger.
        \textbf{D} Histogram of the quadrature weights shown in panel \textbf{C}.
    }
    \label{fig:surface_convergence}
\end{figure}

\subsection{Convergence on a Torus} \label{ssec:surface_convergence}
We next consider convergence results for quadrature and neural field simulations on a torus. For these experiments, we use spiral nodes for the quadrature (See \Cref{sec:numerical_tests}). 
While the algorithm does not require such regularity, using spiral nodes simplifies the setup and reduces error variability.

To test the surface RBF-QF algorithm directly, we use the function $f(x, y, z) = \sin(7x) + 1$ on the torus $D_T$, with $\phi(r) = r^3$, stencil size $k=21$, and an a polynomial basis appended up to a specified degree. 
While we lack a direct mesh radius on the surface, the nodes are well spaced, so we
use $n^{-1/2} = \mathcal{O}(h)$ as a proxy measure of resolution~\cite{REEGER2016,
REEGER2018}. \Cref{fig:surface_convergence}{A} shows that for all polynomial
degrees tested, convergence exceeds the expected rate
$\mathcal{O}(n^{-\text{deg}/2})$.

To test our method of lines neural field simulation on a surface, we construct a
manufactured solution analogous to that used in the flat case. Specifically, we let
\begin{align*}
    u(t, \vecx) &= f^{-1}\left[ \text{Gauss}(\psi(\vecx), \vecx_0; \sigma=11/10) + 1/10\right],
        \qquad 
        \vecx_0 = [\cos(t)/5, \sin(t)/5]^T
\end{align*}
solve \cref{eqn:manufactured_system} with $f$ as in \cref{eqn:firing_rate} and weight kernel
\[
    w(\vecx, \vecy) = 
    \text{Gauss}(\psi(\vecx), \psi(\vecy); \sigma=1/40) \frac{1}{(R +
    r(\vecy)\cos \theta(\vecy))}.
\]
This mirrors the previous solution but is defined on the surface domain $D_{2\pi}$ (in $\phi$-$\theta$ coordinates), with the weight kernel adjusted by the inverse Jacobian of the torus parametrization ($[R+r \cos \theta]^{-1}$). Although the analytical solution is identical, the quadrature weights differ due to the use of surface-specific integration techniques.

Results are shown in \Cref{fig:surface_convergence}. 
Panel {A} shows quadrature convergence on the torus using the test function
$f(x, y, z) = \sin(7x)+1$ with relative error decreasing as
$n^{-1/2}$ and higher decay rates for higher-degree appended polynomials.
Panel~{B} reports convergence for the manufactured neural field solution,
again showing high-order accuracy.
Panel~{C} visualizes the quadrature weights across the torus surface for $n =
32000$, while {D} presents a histogram of these weights.
Weights are constant along lines of constant $\theta$ due to the symmetry of the node
placement; their distribution reflects both curvature effects and histogram binning.

\section{Showcase of neural field dynamics on curved geometries}\label{sec:showcase}
Thus far, we have described and tested a numerical method for solving neural field
equations on smooth, closed manifolds. In this section, we demonstrate its utility
with simulations that reveal rich spatiotemporal dynamics on curved geometries.
Each example extends classical results from planar or one-dimensional settings to
non-Euclidean domains, identifying impacts and considerations of curvature on
solution dynamics, opening new avenues for mathematical and scientific investigation.
All computations use the same geodesic distance approximation, described in the
\hyperref[sec:code]{\textbf{Code Availability}} section.

\subsection{Labyrinthine patterns on a deformed sphere}\label{ssec:blood_cell}
Previously, Coombes et al~\cite{COOMBES2012} investigated a neural field model in two
spatial dimensions that exhibits {\em labyrinthine patterns} under appropriate
initial condition. 
These patterns emerge when a lateral inhibitory weight kernel is used with an
unstable, near symmetric bump initial activity profile. Owing to the $D_4$ symmetry
of the dominant instability, four arms of a cross extend outward from the bump and
begin to branch, producing a single connected region of high neural activity, forming
a complex labyrinth of thin, repeatedly branching corridors.

We perform a similar numerical experiment on a family of deformed spheres,
demonstrating that surface curvature has a significant effect on the qualitative
properties of the resulting labyrinthine patterns.  
Unlike the original experiment by Coombes et al.~\cite{COOMBES2012}, which employed a
discontinuous Heaviside firing-rate function and a non-smooth weight kernel, we use
smooth, qualitatively similar functions that still give rise to complex
spatiotemporal activity. Specifically, we define the weight kernel and firing-rate
function as
\begin{align}
w(\vecx, \vecy) &= A_e \, \text{Gauss}(\vecx, \vecy, \sigma_e) - A_i\,
\text{Gauss}(\vecx, \vecy, \sigma_i), \label{eq:winh} \\
f(u) &= 
\begin{cases}
0, & u < 0.06, \\
p(u), & 0.06 \le u < 0.54, \\
1, & u \ge 0.54,
\end{cases} \label{eq:smoothf}
\end{align}
where $A_e=5$, $A_i = 5$ $\sigma_e = 0.05$, $\sigma_i = 0.1$. Although $A_e=A_i$, the
kernel satisfies $w(\mathbf{x}, \mathbf{x})>0$ since the Gaussian is normalized --
its peak is higher for smaller $\sigma$ -- producing the desired lateral inhibitory
property. The function $p(u)$ is a $9$th-order polynomial chosen to ensure that \( f \)
is four times continuously differentiable, yielding a smooth spline approximation to
the Heaviside function \( H(u - 0.3) \), remaining constant outside a narrow
transition interval centered at the threshold.

\begin{figure}
    \centering
     \href{https://figshare.com/articles/media/Animations/28791911?file=53662766}{\includegraphics[width=1.0\textwidth]{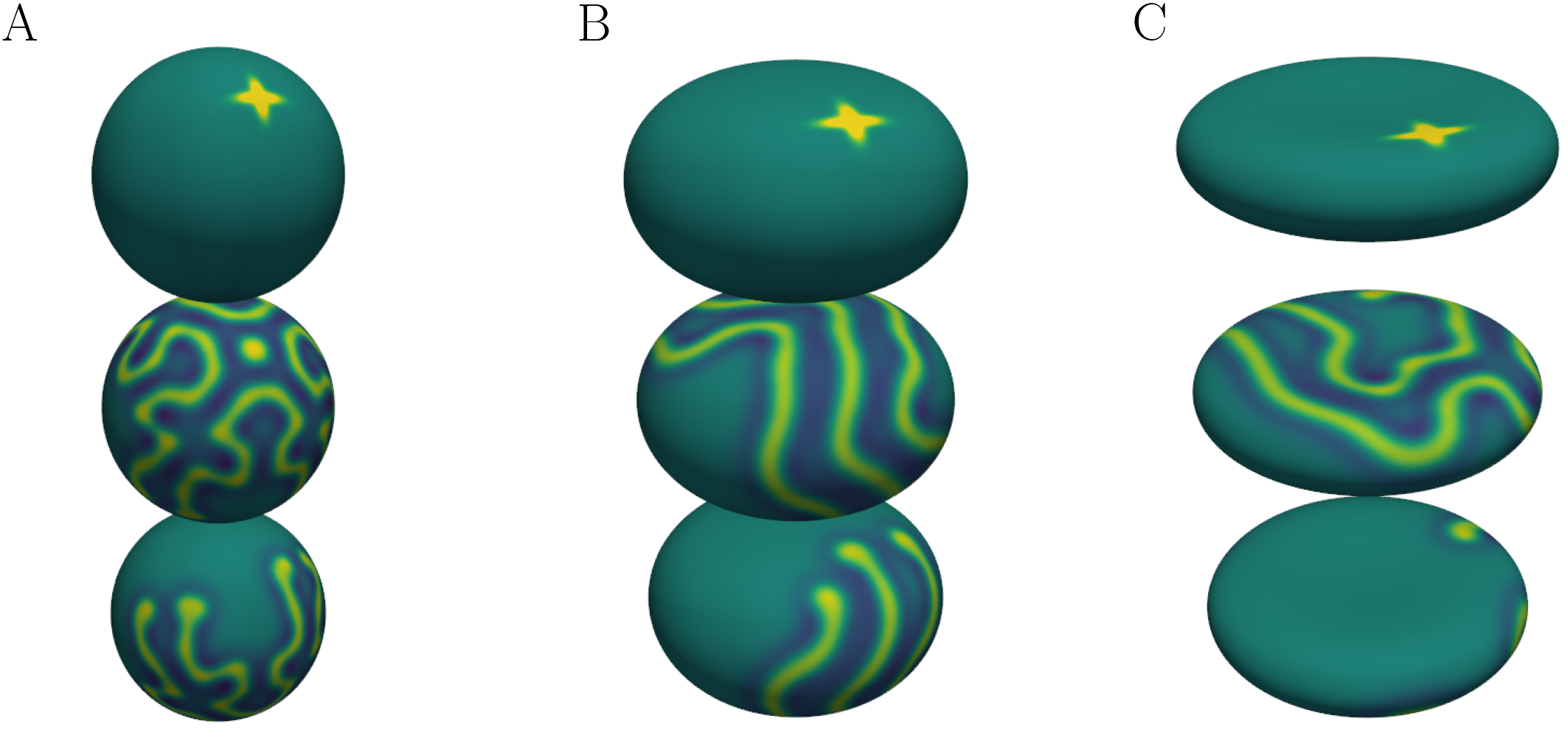}}
    \caption{
        \textbf{Effect of Surface Deformation on Labyrinthine Patterns.}  
Labyrinthine patterns on a deformed spherical surface \( D_\gamma \), with increasing deformation parameter from left to right:  
\textbf{A.} \( \gamma = 0 \) (perfect sphere),  
\textbf{B.} \( \gamma = 0.4 \) (moderate deformation),  
\textbf{C.} \( \gamma = 0.8 \) (flattened, blood-cell-like shape).  
An initial unstable, cross-shaped activity pattern (top) evolves to meandering spatiotemporal patterns (middle: north pole viewpoint; bottom: south pole viewpoint). 
Yellow indicates high neural activity (\( u > 0 \)), green denotes baseline activity (\( u \approx 0 \)), and dark green/blue indicates suppressed activity (\( u < 0 \)).  
Surface curvature has a profound effect on the form and complexity of the resulting patterns. Click image to view full animation 
(\href{https://figshare.com/articles/media/Animations/28791911?file=53662769}{Movie S1A}, 
\href{https://figshare.com/articles/media/Animations/28791911?file=53662775}{Movie S1B}, 
\href{https://figshare.com/articles/media/Animations/28791911?file=53662766}{Movie S1C}).
}
    \label{fig:blood_cell}
\end{figure}

The underlying surface \( D_\gamma \) is a one-parameter family of deformed spheres,
with deformation governed by \( \gamma \in [0, 1) \). When \( \gamma = 0 \), the
surface is a standard sphere; as \( \gamma \) increases, the geometry is increasingly
compressed along the vertical axis, reducing the pole-to-pole Euclidean distance to
\( 2(1 - \gamma) \). The surface is implicitly defined by
\[
1 = x^2 + y^2 + z^2 \left(1 - \frac{\gamma}{1 + (x^2 + y^2)/2.89} \right)^{-1},
\]
and is visualized in \Cref{fig:blood_cell} for \( \gamma = 0, 0.4, 0.8 \).

We initialize the activity \( u(0,\vecx) \) on the undeformed sphere using the function
\[
u(0,\vecx) = 5 \cdot \exp \left( -10 \left[ \left(\cos\left(4 \arctan\left(\frac{y}{x}\right)\right) + 3 \right) \sqrt{x^2 + y^2} \right]^2 \right) H(z),
\]
which defines a symmetric, cross-shaped pattern localized in the northern hemisphere.  
To orient this pattern appropriately, we apply a rigid rotation that maps the north
pole \([0, 0, 1]\) to the point \([0.5, 0.3, \sqrt{0.5^2 + 0.3^2}]\), aligning the
initial condition with the desired surface region before the deformation is applied.

We then simulate the neural field on each deformed surface using the kernel basis
function \( \phi(r) = r^3 \), a stencil size of \( k = 32 \), and append third-degree
polynomials for accurate quadrature. Time integration is performed using the
Adams–Bashforth 5 method with a fixed time step of \( \Delta t = 10^{-2} \). The
resulting dynamics, corresponding to various values of the surface deformation
parameter \( \gamma = 0, 0.4, 0.8 \), are shown in \Cref{fig:blood_cell}. The top row
depicts the initial condition, rendered from a viewpoint where the {\em north pole}
is visible.  
The middle row shows the activity at time \( t = 200 \) from the same viewpoint,
while the bottom row presents the state at \( t = 200 \) with the surface inverted
through the \( x \)-\( y \) plane to reveal the {\em south pole}.

Labyrinthine corridors remain connected but their wandering termini tend to veer away
from regions of high curvature, illustrate how surface geometry profoundly influences
neural field dynamics. This is akin to the pinning of propagating
waves~\cite{bressloff2001traveling,coombes2011pulsating} or attraction/repulsion of
localized activity~\cite{kilpatrick2013wandering} observed by introducing
inhomogeneities into the weight kernel of neural fields on flat domains.

\subsection{Traveling spot steered by surface bumps}

Our next example demonstrates how surface curvature influences the trajectory of
traveling spot solutions in a neural field model with synaptic
depression~\cite{bressloff2011two,shaw2024representing}. The model consists of
two coupled equations: one for the neural activity $u$, and one for the synaptic
efficacy $q$, ranging from 0 (depleted) to 1 (fully available). The dynamics are
given by
\begin{align}
    \partial_t u &= -u + \int_D w(\blank, \vecy) q(\blank, \vecy) f[u(\blank, \vecy)] \ d\vecy,
    \qquad \tau \partial_t q = 1 - q - \beta q f[u]. \label{eq:nfsyndep}
\end{align}
We use the laterally inhibitory weight kernel \cref{eq:winh} with $A_e = 5$, $A_i =
7$, $\sigma_e = 0.05$, $\sigma_i = 0.1$. Such weights often produce {\em spot}
solutions (localized circular active regions, also called pulses or bumps) in planar
neural fields without adaptation ($\beta = 0$).

\begin{figure}
    \centering
    \href{https://figshare.com/articles/media/Animations/28791911?file=53662781}{\includegraphics[width=1.0\textwidth]{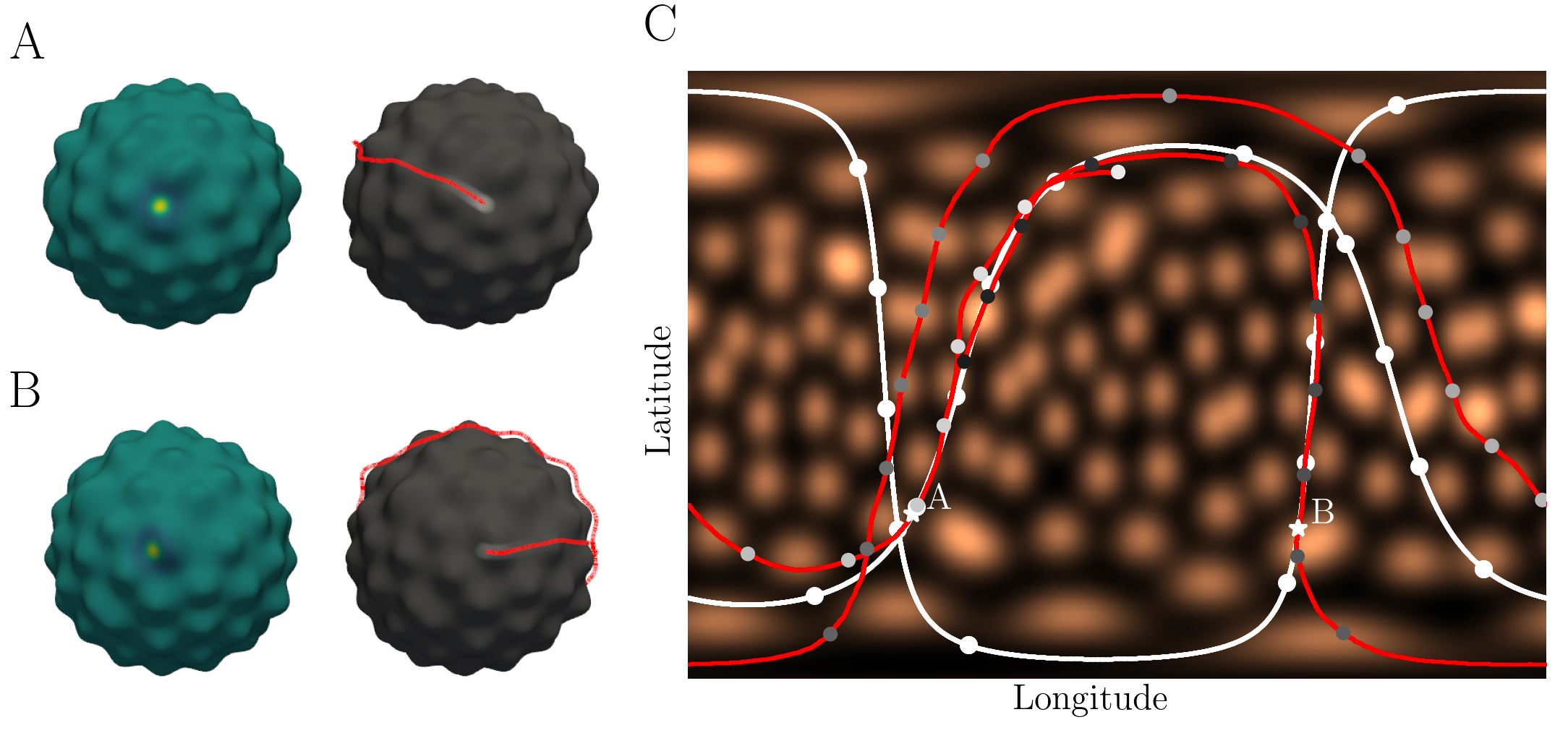}}
        \caption{\textbf{Traveling spot steered by bumps in curvature.}  
\textbf{A–B:} Snapshots from a simulation of a traveling spot on a bumpy sphere. Each shows the activity variable (left) and synaptic efficacy (right). The trajectory of the spot up to that
point is shown as a red curve hovering above the surface.
\textbf{C:} Projection of the bumpy surface into spherical coordinates, with
longitude on the horizontal axis, latitude on the vertical axis, and surface
elevation encoded by heatmap. Traveling spot trajectory (red) sometimes hugs but also
deviates from great-circle tangents (white) aligned spot position in panels
\textbf{A} and \textbf{B}, showing how bumps deflect motion. Click image to view full
animation
(\href{https://figshare.com/articles/media/Animations/28791911?file=53662781}{Movie
S2}).
}
    \label{fig:bumpy_sphere}
\end{figure}

Incorporating synaptic depression ($\beta > 0$), as with other forms of adaptation~\cite{owen2007bumps}, causes spots to travel. Introducing depression ($q < 1$) on one side of the spot, effective lateral inhibition is asymmetric and activity will increases on the side of the spot farthest from the depressed region~\cite{bressloff2011two,shaw2024representing}. The spot then propagates, leaving a trail of synaptic depression in its wake. 

We probe how surface curvature affects trajectories of such traveling spots by
considering a surface that we refer to as a {\em bumpy sphere}
(\Cref{fig:bumpy_sphere}), defined by the equation:
\[
    1 + \sum\limits_{\vecx_i} \left(\frac{2\pi}{10^2}\text{Gauss}(\vecx, \vecx_i, 1/10)\right) \tfrac{1}{10} = x^2 + y^2 + z^2
\]
where $\{\vecx_i\}_{i=1}^{100}$ are a set of 100 {\em bump centers} (randomly chosen,
though roughly evenly spaced) on the unit sphere (See remarks in
\hyperref[sec:code]{\textbf{Code Availability}} section).

\Cref{fig:bumpy_sphere}{A} and {B} each show two snapshots from the
same simulation, each depicting two views of the bumpy sphere. The left view is
colored by the activity variable, while the right shows synaptic efficacy in
grayscale. Both are rotated to center the spot in view, and the grayscale surface
includes a red curve traving the trajectory of the spot up to that time.
\Cref{fig:bumpy_sphere}{C} presents the surface in spherical coordinates, with latitude and longitude along the vertical and horizontal
axis, and color representing radial elevation. The red curve shows the full
trajectory of the traveling spot. If the surface were perfectly spherical, the spot
would follow a straight path along a great circle, akin to what was found on the
planar case \cite{bressloff2011two}. Instead, we observe consistent
deviations in the trajectory that arise from geometric inhomogeneities. 
Two white curves show tangent great circles at times shown in
panels {A} and {B}, illustrating how the trajectory is deflected by geometric
inhomogeneities. This reenforces the point that surface irregularities creates a
similar potential surface which shape the dynamics of evolving spatiotemporal
solutions, akin to those created in neural fields on flat domains with weight
inhomogeneities~\cite{poll2015stochastic,kuehn2019gradient}.

\subsection{Labyrinthine patterns on a human cortex}

\begin{figure}
    \centering
    \href{https://figshare.com/articles/media/Animations/28791911?file=53662760}{\includegraphics[width=1.0\textwidth]{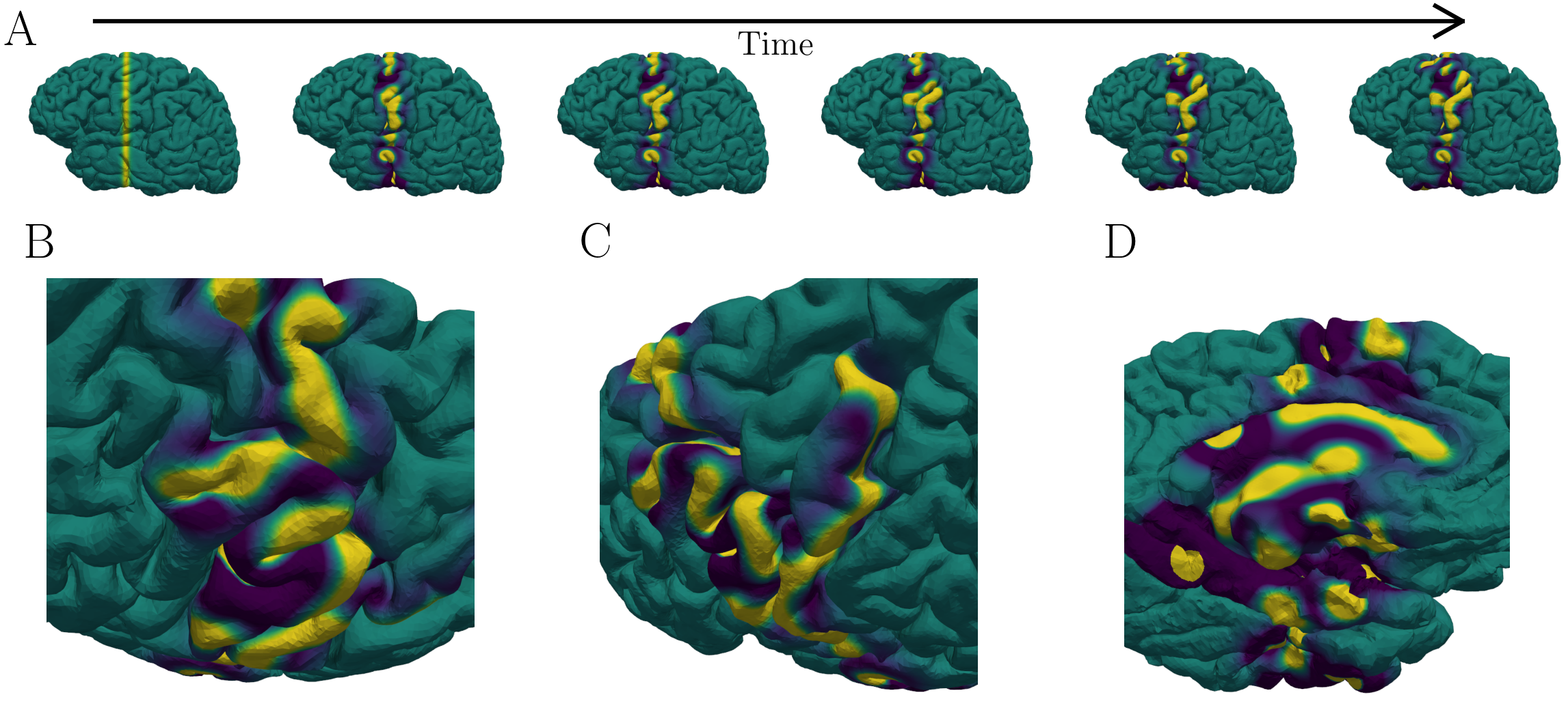}}
    \caption{{\bf Labyrinthine dynamics on a cortex.} Here we show a simulation on
      the left hemisphere of a realistic human cortex using a laterally inhibitory
      weight kernel that favors the formation of labyrinthine patterns.
\textbf{A}: Snapshots at evenly spaced early time points show the initial stripe of
activity breaking into multiple labyrinth-like corridors that begin to grow and
evolve.
\textbf{B–D}: Final-time snapshots from different viewing angles highlight distinct
cortical regions: \textbf{B.}~frontal lobe, \textbf{C.}~parietal lobe, and
\textbf{D.}~limbic cortex/medial surface. Click image to view full animation
(\href{https://figshare.com/articles/media/Animations/28791911?file=53662760}{Movie S3}).
    }
    \label{fig:cortex}
\end{figure}

Finally, we build a simulation showcasing the evolution of a neural field on a
cortical surface extracted from human data (\Cref{fig:cortex}). See
\hyperref[sec:code]{\textbf{Code Availability}} for a link to the repository
containing this mesh, representing the left hemisphere of a human cortex. It was
generated using MNE-Python, which integrates FreeSurfer’s anatomical reconstruction
pipeline to segment T1-weighted MRI scans and extract detailed cortical surfaces. The
resulting triangulated mesh captures the geometry of the pial surface, including
sulci and gyri, enabling anatomically realistic neural field simulations. We use the
same laterally inhibitory weight function defined in \cref{eq:winh} with $A_e =
A_i = 5$, $\sigma_e = 3$, $\sigma_i = 6$, along with the geodesic distance metric and
parameters from \Cref{ssec:blood_cell}. The initial condition consists of a band of
activity, which rapidly fragments into localized regions of varying size shaped by
the surface curvature. The resulting dynamics include both stationary spots and
labyrinthine corridors that tend to follow gyri and sulci—regions aligned with
locally minimal curvature.

The mesh used here is visually compelling and
anatomically detailed, though originally intended for visualization rather than
numerical simulation. It exhibits some geometric irregularities -- including a small
duplicated region, uneven node density, and inconsistencies in surface normals. This
leads to rapidly oscillating quadrature weights, including some large negative
values.
\rev{
However, this is not indicative of an inherent instability of the RBF quadrature
method. On other domains, including the flat-domain tests in
Section~\ref{ssec:square} and the torus in Section~\ref{ssec:surface_convergence},
the method produces stable, accurate results with only small negative weights, even
over long simulations. These observations demonstrate that the method is robust when
applied to high-quality, well-distributed node sets. The oscillatory weights in the
cortex case stem from the severe irregularity of the publicly available mesh we
employed, rather than from the quadrature scheme itself.
}
Nonetheless, the resulting dynamics reveal rich qualitative
structure that underscores the potential of this framework and motivates further
refinement and study.

Although a rigorous analytical treatment of these effects on curved surfaces remains
an open question, these results highlight the potential for geometry to shape the
dynamics of cortical activity. Our framework provides a powerful computational tool
for probing such curvature-driven effects in structured neural field models.

\section{Conclusion}

We have presented a high-order, mesh-flexible solver for neural fields on smooth, closed surfaces using RBF-based interpolation and quadrature. The method is numerically stable, accurate, and requires only a triangulated mesh and approximate vertex normals. Unlike spectral methods, which rely on structured domains, our approach—like finite element methods—supports arbitrary geometries without requiring element construction. Although the resulting quadrature matrices are sparse, pairwise interactions lead to $\mathcal{O}(n^2)$ complexity in the number of nodes. This scaling motivates the use of high-order schemes that achieve accuracy with fewer degrees of freedom, especially in cortex-scale modeling or inverse problems involving kernel learning.

Simulations on bumpy spheres and realistic cortical surfaces show how geometry can steer and constrain activity, extending phenomena observed in planar neural fields with inhomogeneous coupling (e.g., wave slowing, deflection, or pinning) to non-Euclidean domains. In \Cref{fig:blood_cell}, a gyrus-like ridge steers a labyrinthine wave pattern along low-curvature paths until repulsion from adjacent corridors forces a transition; similar behavior is observed on real cortical gyri in \Cref{fig:cortex}. The bumpy sphere simulation in \Cref{fig:bumpy_sphere} further demonstrates that curvature can deflect traveling spot trajectories. While its effects on wave speed and stability remain unclear, prior work suggests that curvature can pin or disrupt waves~\cite{bressloff2001traveling}. Multi-spot simulations reveal curvature-modulated crowding and spot annihilation~\cite{kilpatrick2013wandering}, raising broader questions about how surface geometry and the excitatory–inhibitory balance of the kernel interact to steer dynamics, and whether a critical angle of incidence maximizes deflection. These findings motivate future reductions to effective equations and further analysis of curvature-driven pinning, transitions, and stability~\cite{bressloff2019stochastic,martin2018numerical}.


\rev{
Surface differential operators (e.g., diffusion or advection) commonly arise in
neural field models with local dynamics~\cite{kneer2014nucleation,BASPINAR2023}.
RBF-based finite difference methods provide high-order, geometry-flexible
approximations of such
operators~\cite{petras2018rbf,shankar2015radial,lehto2017radial,shaw2019radial,piret2012orthogonal,fuselier2013high},
with the same system matrix used for quadrature weights also yielding finite
difference weights—offering computational savings when coupling local and nonlocal
dynamics. Other meshfree approaches, such as partition of unity
methods~\cite{babuska1997partition,wendland2002fast} and moving least
squares~\cite{levin1998mls,fries2010mls}, share similar flexibility for irregular
node layouts and could be adapted to neural field models on surfaces. We focus on
RBF-based methods for their direct unification of interpolation, quadrature, and
differentiation, noting that recent RBF-FD advances in stabilization and adaptive
refinement~\cite{BAYONA2017,FLYER2016} further expand the toolkit for high-order,
geometry-flexible PDE solvers.
}

Our method depends only on surface geometry and supports arbitrary kernels, offering
a promising platform for pairing with experimental data to infer connectivity. This
opens the door to data-driven modeling, model inversion, and further theoretical
exploration of how cortical geometry shapes large-scale neural activity.

\subsection*{Code Availability} \label{sec:code}
The code used to generate the numerical simulations and figures can be found in the repository \href{https://github.com/shawsa/neural-fields-rbf}{\texttt{www.github.com/shawsa/neural-fields-rbf}}.
For curved domains, the first order
approximation to the geodesic distance via the Fast Marching
Algorithm~\cite{sethian1996fast} can be found in the MeshLib library:
\href{https://meshlib.io/}{\texttt{meshlib.io}}. The realistic cortical mesh was adapted from the MNE python library:
\href{https://mne.tools/stable/index.html}{\texttt{mne.tools}}.

\subsection*{Supplementary Material}

\begin{figure}
  \centering
  \includegraphics[width=1.0\textwidth]{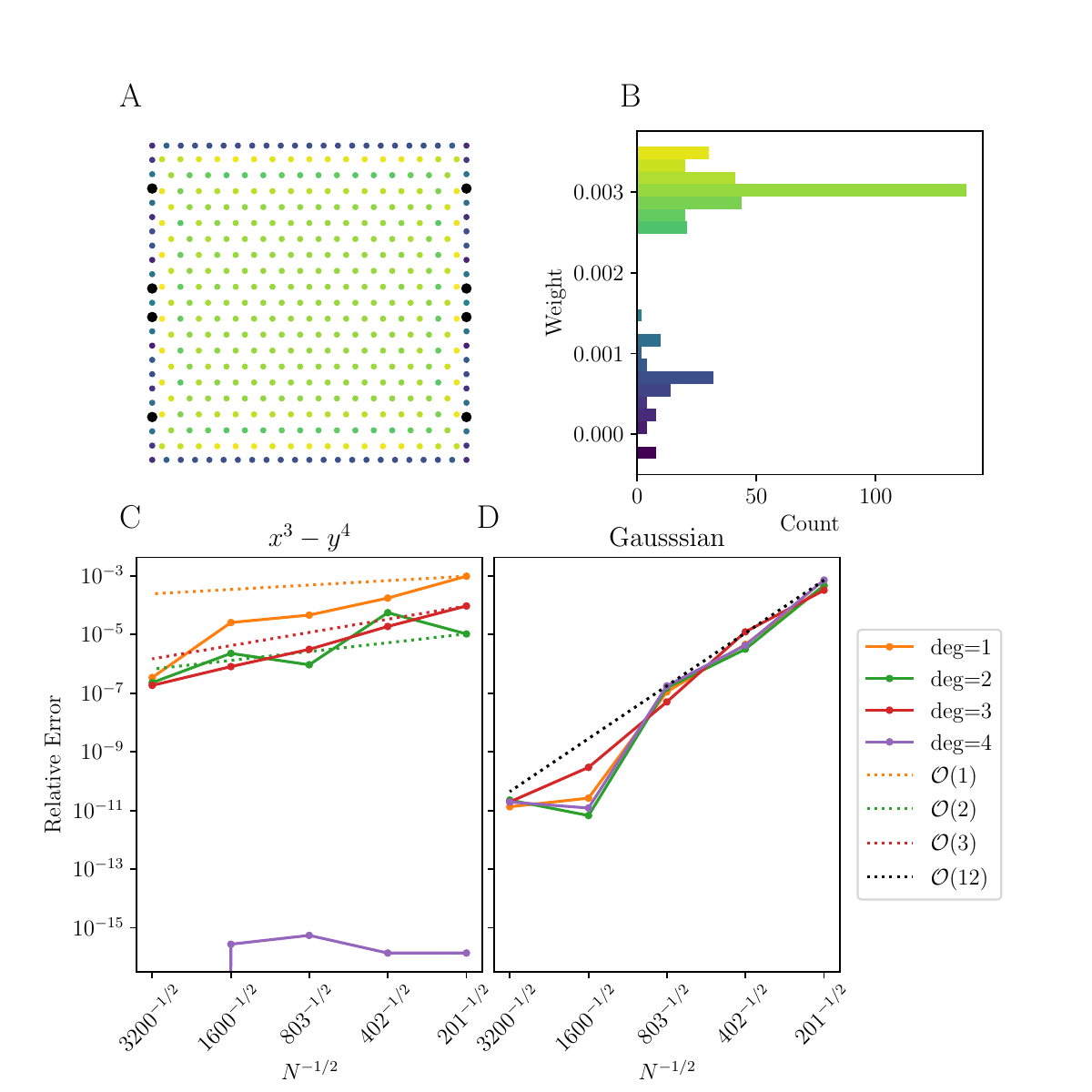}
  \caption{
  Quadrature on the unit square using a regular triangular grid of collocation nodes, the RBF $\phi(r) = r^3$, and stencil size $k=30$.
  \textbf{A.} Quadrature nodes colored by weight (negative weights are shown as larger black dots).
  \textbf{B.} Histogram of quadrature weights.
  \textbf{C.} Convergence of the quadrature rule for a degree-4 polynomial $x^3-y^4$.
  \textbf{D.} Convergence for a Gaussian $f(x,y)=\exp(-10[(x-\tfrac{1}{2})^2+(y-\tfrac{1}{2})^2])$. Relative error is plotted against $N^{-1/2}$, corresponding to the node length scale.
  }
  \label{fig:hex_quad}
\end{figure}

We present convergence results on the unit square, where collocation nodes are placed
on a regular triangular grid in the interior and equally spaced along the boundary.
Figure~\ref{fig:hex_quad}A,B shows the node placement and corresponding quadrature weights.
Convergence tests are reported in Figure~\ref{fig:hex_quad}C,D for two representative
functions: a degree-4 polynomial and a Gaussian. For the polynomial, the quadrature
is exact when augmenting the basis with degree-4 polynomials, as expected.
The Gaussian is given by
\[
    f(x, y) = \exp\bigg(-10\big[(x-\tfrac{1}{2})^2 + (y-\tfrac{1}{2})^2\big]\bigg)
\]
decays rapidly to zero near the boundary. Because this function is smooth and effectively
supported away from the edges, the quadrature achieves spectral convergence.
This behavior is directly analogous to the classical result that the trapezoidal rule
attains spectral accuracy for smooth periodic functions~\cite{trefethen2000spectral}.


\section{Convergence on other geometries}
We next present convergence results for quadrature on curved surfaces, specifically
the unit sphere and a Dupin cyclide. 
For the sphere, collocation nodes are chosen from icosahedral-based point sets \cite{teanby2006icosahedron, hardin2016comparison}, which provide nearly uniform coverage.
For the cyclide, we consider the ring case with parameters
\[
a = 1,\quad b = 0.98,\quad c = 0.1983,\quad d = 0.5,
\]
and construct a non-random triangular mesh in the $(\phi,\theta)$ parameter space. 

\begin{figure}[t!]
  \centering
  \includegraphics[width=1.0\textwidth]{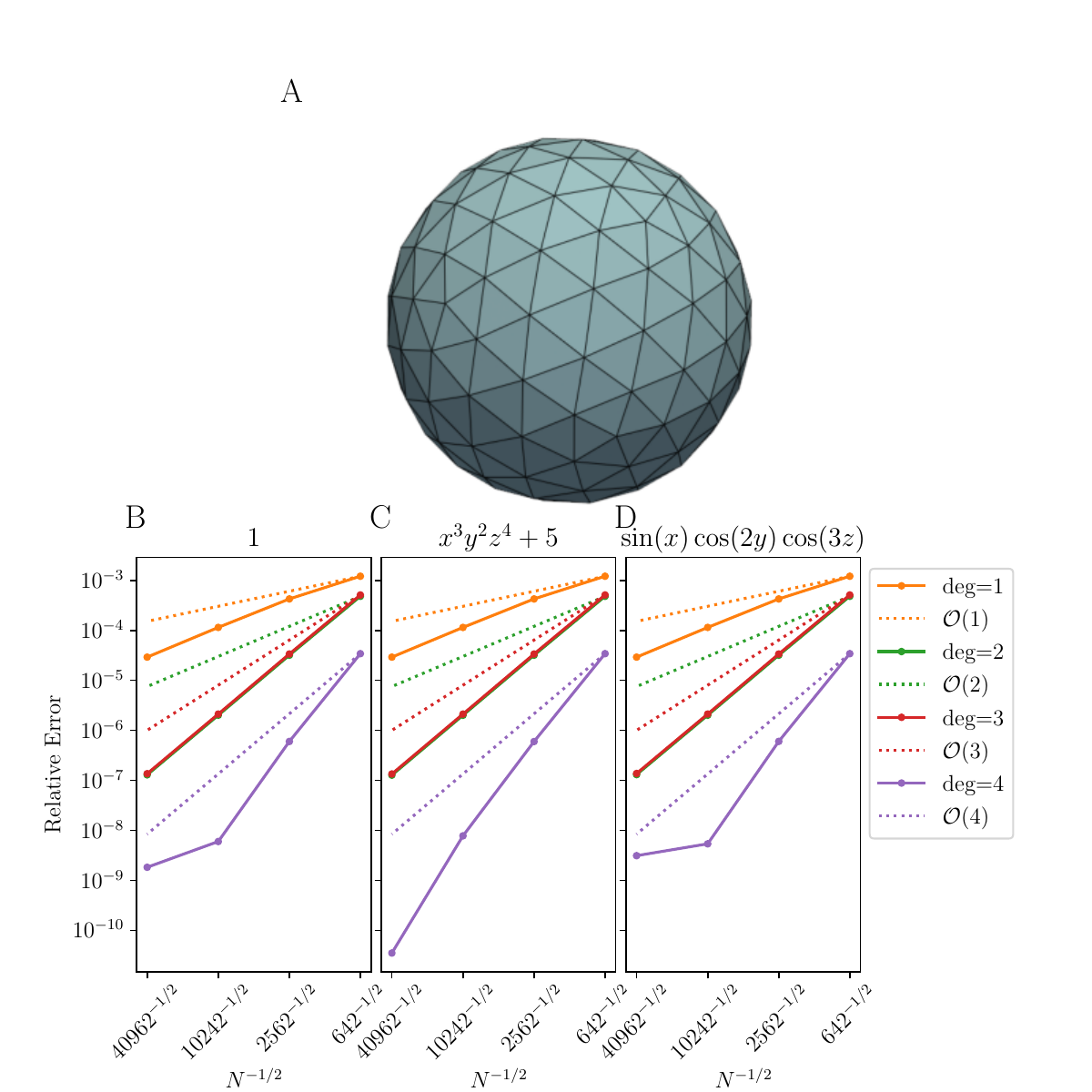}
  \caption{
  Convergence of quadrature on the unit sphere using icosahedral nodes.
  \textbf{A.} Example of a triangular mesh from icosahedral nodes placement.
  \textbf{B-D.} Relative quadrature error vs $N^{-1/2}$ (node length scale) for three test functions: a constant (1), a polynomial $x^3y^2z^4 + 5$, and a trigonometric function  $\sin (x) \cos(2y) \cos (3z)$.
  }
  \label{fig:sphere_quad}
\end{figure}

\begin{figure}[t!]
  \centering
  \includegraphics[width=1.0\textwidth]{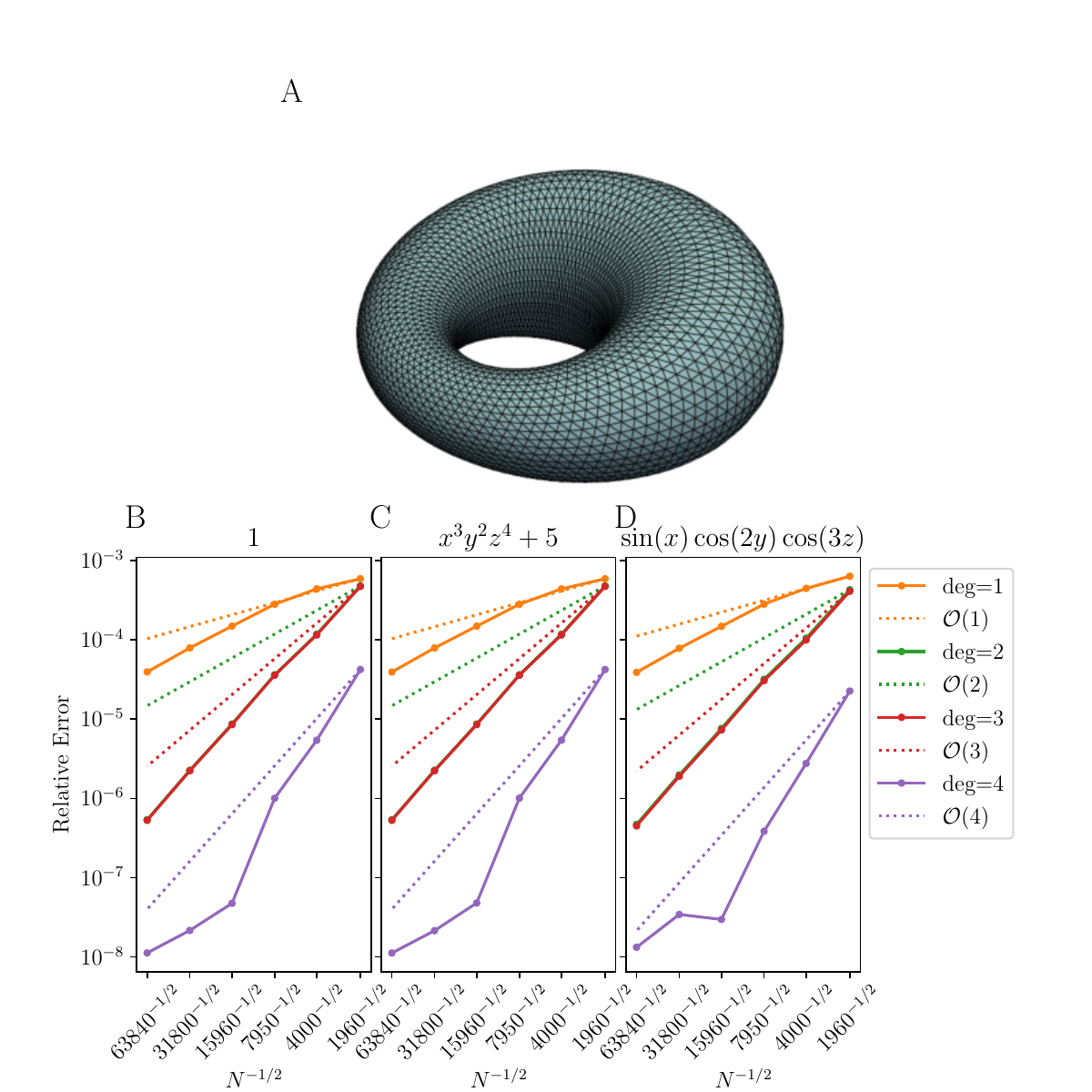}
  \caption{
  Convergence of quadrature on a ring cyclide with parameters $a = 1$, $b = 0.98$, $c = 0.1983$, $d = 0.5$, using a non-random triangular mesh in $(\phi, \theta)$ parameter space.
  \textbf{A.} Example of the cyclide mesh.
  \textbf{B-D.} Relative quadrature error vs $N^{-1/2}$ (node length scale) for three test functions: a constant (1), a polynomial $x^3y^2z^4 + 5$, and a trigonometric function  $\sin (x) \cos(2y) \cos (3z)$.
  }
  \label{fig:cyclide_quad}
\end{figure}

Figures \ref{fig:sphere_quad}A and \ref{fig:cyclide_quad}A illustrate representative meshes. 
Panels B--D in each figure report relative quadrature error versus $N^{-1/2}$ for three test functions: a constant $f(x,y,z)=1$, a polynomial $f(x,y,z)=x^3y^2z^4+5$, and a trigonometric function $f(x,y,z)=\sin(x)\cos(2y)\cos(3z)$. 
As predicted, we observe convergence at least has high as the degree of appended polynomial for all three test functions.

\bibliographystyle{siamplain}
\bibliography{references_abbreviated}

\end{document}